\documentclass[11pt,draft,twoside,a4paper]{book}
\usepackage{amsmath, amssymb, amsthm, amsfonts}
\usepackage{indentfirst}
\usepackage[iso-8859-7]{inputenc}
\usepackage{amsmath}
\usepackage{amsfonts}
\usepackage{latexsym}
\usepackage{graphics}
\usepackage{graphicx}
\newcommand{\tl }{\textlatin }

\begin{document}
\small
\hsize=14.5truecm                 
\newtheorem{definition}{Definition}[chapter]
\newtheorem{lemma}{Lemma}[chapter]
\newtheorem{theorem}{Theorem}[chapter]
\newtheorem{proposition}{Proposition}[chapter]
\newtheorem{corollary}{Corollary}[chapter]
\newtheorem{remark}{Remark}[chapter]%
\newcommand{\pa}{\hskip 0.6truecm}
\newcommand{\HRule}{\rule{\linewidth}{1mm}}
\newcommand{\fr }{\tl{Fourier}\;}
\newcommand {\cpi}{\frac{1}{(\sqrt{2\pi})^n}}
\newcommand{\abs}[1]{\left\vert#1\right\vert}
\def\e{\epsilon}
\def\C{\mathbb{C}}
\def\Q{\mathbb{Q}}
\def\N{\mathbb{N}} 
\def\Z{\mathbb{Z}} 
\def\R{\mathbb{R}} 
\def\Rn{$\R^n\;$}
\def\ApO{$A^p(\Omega)$ }
\def\Ap {$A^p(\Delta)$ }
\def\f{$f$}
\def\limit{\hbox{\space \raise -.2mm\hbox{$ \longrightarrow \atop {n\rightarrow \infty}$} \space}}

\setlength{\oddsidemargin}{0.20in}
\setlength{\evensidemargin}{0.20in}
\setlength{\topmargin}{0.0in}

\thispagestyle{empty}
\medskip
\vspace*{\stretch{2}}
\medskip

\centerline {\huge\textbf{Bergman Spaces Seminar}}
\bigskip
\medskip

\vspace*{\stretch{2}}
\centerline  {\Large\textsc{Lecture Notes}}
\medskip

\vspace*{\stretch{4}}

\centerline {\textsc{ Department of Mathematics}}
\bigskip
\centerline {\textsc{University of Crete}}
\medskip
\centerline {\textsc{Heraklion 2005}}
\vspace*{\stretch{2}}

\thispagestyle{empty}
\newpage
\thispagestyle{empty}

\vspace*{3in}
\tableofcontents

\newpage
\thispagestyle{empty}

\chapter {Introduction to Bergman Spaces }
\label{ch1}
 
\section {Basic Definitions}
\medskip

We begin our discussion with the Bergman Spaces on the unit disc of the complex plane
\begin{eqnarray}\label{eq1.1}
\Delta=\{z\in \C : |z|<1 \} .
\end{eqnarray}
These are defined as 
\begin{eqnarray}\label{eq1.2}
A^p (\Delta)=\{ f\in H(\Delta) : \|f\| _{A^p (\Delta ) }  ^p = \int _{\Delta}|f(z)|^p dm(z)<\infty \}
\end{eqnarray}
where $dm(z)=\frac{1}{\pi}dA(z)$ is the two dimensional Lebesgue measure normalized in $\Delta$ and $H(\Delta)$ is the space of analytic functions on the unit disc.

However, the story begins with Stephan Bergman and around 1970 ("The Kernel function and Conformal Mapping") there's already been some progress in the study of the spaces
\begin{eqnarray}\label{eq1.3}
A^p (\Omega)=\{ f\in H(\Omega) : \|f\|_{A^p((\Omega))} ^p =\int _{\Omega}|f(z)|^p dA(z)<\infty \},\; 0<p<\infty
\end{eqnarray}
where $\Omega \subseteq \C$ is an open connected set, for the case $p=2$. The interest was then focused on the case $\Omega=\Delta$.

The study of the Bergman Spaces is inspired from that of the Hardy Spaces
\begin{eqnarray}\label{eq1.4}
H^p (\Delta)  =  \Big\{ f\in H(\Delta) : \|f\| _{H^p (\Delta ) }  ^p = \sup_{r\in[0,1)}\frac{1}{2\pi}\int _{0} ^{2\pi}|f(re^{i\theta}))|^p d\theta <\infty \Big\}, \; 0<p<\infty  
\end{eqnarray}
and
\begin{eqnarray}\label{eq1.5}
H ^ {\infty} (\Delta )&=& \{ f\in H(\Delta ) : \| f \| _{H^\infty(\Delta )} = \sup _{z \in \Delta} |f(z)|<\infty \} 
\end{eqnarray}
since $H^p(\Delta)\subseteq A^p (\Delta)$. Observe however that $H^\infty(\Delta) = A^\infty(\Delta)$ so hencefotrth we will restrict our interest to the case $p<\infty$.

The study of $H^p$ spaces begun from Hardy between 1915 and 1930 and and then was continued with great interest (around 1960, Lenarnt Carleson and then Shapiro and Shields solved the so called "universal interpolation problems").

As it turned out, the Bergman spaces $A^p(\Delta)$ behave quite differently from the Hardy spaces, and the study of Bergman Spaces remained still until 1990. Then Hedenmalm (in $A^2$) and Duren, Khavinson, Shapiro, Sundberg, Sheip, Aleman (in $A^p$) gave significant results. As a consequence, the study of the spaces $A^p(\Delta)$ was rapidly evolved in the last 15 years.
Lately there is interest in the study of the spaces $A^p(\Omega)$ and of the Bergman spaces on the unit ball.

\medskip


\section {Growth of $A^p$ Functions}
\medskip

We will see in this section a description of the growth of $A^p$ functions and some basic consequences. We will state this result in  the more general context of the spaces $A^p(\Omega)$, where $\Omega \subseteq \C$ is an open connected set.
\begin{proposition}\label{pr1.1}
Let $f\in A^p(\Omega)$, $0<p<\infty$. Then, for every $z\in \Omega$,

\begin{equation}\label{eq1.6}
|f(z)|\leq \frac{1}{\pi ^{\frac{1}{p}}}\delta(z)^{-\frac{2}{p}}\|f\|_{A^p(\Omega)}
\end{equation}
where $\delta(z)=dist(z,\partial \Omega)$.
\end{proposition}

\begin{proof}
Fix $z \in \Omega$ and set $\delta=\delta(z)$. We consider the disc $\Delta(z,\zeta) =\{\zeta \in \Omega : |\zeta - z|<\delta\}$. Since $f$ is subharmonic in $\Omega$ we have that
$$|f(z)|^p \leq  \frac{1}{2\pi} \int _0 ^{2\pi} |f(z+re^{i\theta})|^p d \theta , \;\; 0\leq r <\delta $$
and so 
$$\int _0 ^\delta |f(z)|^p \, r dr \leq \int _0 ^\delta  \frac{1}{2\pi} \int _0 ^{2\pi} |f(z+re^{i\theta})|^p d \theta \, r dr
.$$ Therefore, 
$$\delta ^2 |f(z)|^p \leq \frac{1}{\pi} \int _{\Delta(z,\zeta)}|f(\zeta)|^p dA(\zeta)\leq \frac{1}{\pi} \int _{\Omega}|f(\zeta)|^p dA(\zeta)$$
which gives
$$ |f(z)|\leq \pi^{-\frac {1}{p}} \delta ^{\frac{-2}{p}} \|f\|_{A^p (\Omega)}$$
which is the desired result.

In the special case of the unit disc $\Delta$ the statement of the theorem becomes 
\begin{eqnarray}\label{eq1.7}
|f(z)|\leq \frac{1}{(1-|z|)^\frac{2}{p}} \|f\|_{A^p(\Delta)}.
\end{eqnarray}
\end{proof}

The analogous description of the growth of the derivatives of a function in $A^p(\Delta)$ is contained in the following Proposition.

\begin{proposition} \label{pr1.2}Let $n$ be a positive integer greater than one and $f\in A^p(\Delta)$, $0<p<\infty $. Then 
\begin{equation}\label{eq1.8}
|f^{(n)}(z)| \leq  \frac{n!2^n 2^\frac{2}{p}}{(1-|z|)^{n+1+\frac{2}{p}}} \|f\|_{A^p(\Delta)}
\end{equation}
\end{proposition}
\begin{proof} Let $r<1$ and set $C=\{\zeta \in \C : |\zeta|=\frac{1+r}{2}\}$. Using Cauchy's integral formula for the derivatives of the function $f$ and the circle $C$ we have that
\begin{eqnarray*}
|f^{(n)}(z)|&=& \Bigg | \frac{n!}{2\pi i} \int _C \frac{f(\zeta)}{\zeta - z}^{n+1} d \zeta \Bigg | \leq 
\frac{n!}{2\pi} \int _ C \frac{|f(\zeta)|}{|\zeta - z|^{n+1}}|d \zeta |  \leq
\frac{n!}{2\pi} \int _ C \frac{|f(\zeta)|}{\big | |\zeta| - |z|\big|  ^{n+1}}|d \zeta |   \\
&\leq & \frac{n!2^{n+1} }{2\pi (1-r)^{n+1}}\int _C |f(z)||d\zeta| \leq 
\frac{n!2^{n+1} 2^\frac{2}{p}}{2\pi (1-r)^{n+1+\frac{2}{p}}} \| f \|_ {A^p(\Delta)} 2 \pi \big( \frac{1+r}{2}\big) \\ 
&\leq & \frac{n!2^n 2 ^\frac{2}{p}}{(1-r)^{n+1+\frac{2}{p}}} \|f\|_{A^p(\Delta)} .
\end{eqnarray*}
Seting $r=|z|$ we get the desired estimate.
\end{proof}

Let us take a look at the consequences implied by the growth of $A^p$ functions described in Proposition \ref{pr1.1}.

First of all it is easy to see that convergence in $A^p(\Omega)$ implies uniform convergence on the compact sets of $\Omega$. While the proof is obvious, we will state this result as a lemma for future reference.

\begin{lemma}\label{lm1.1} Let $\{f_n\}_n$ be a sequence of functions in $A^p(\Omega)$ and $f\in A^p(\Omega)$. Suppose that $f_n$ converges to $f$ in $A^p(\Omega)$. Then $f_n$ converges to $f$
uniformly on the compact subsets of $\Omega $. As a result, $f_n$ converges to $f$ almost everywhere in $\Omega$.
\end{lemma}

We are now ready to show that the spaces $A^p(\Omega)$, $0<p<\infty$, are complete. Of course, $\| \cdot \|_{A^p(\Omega)}$ is a norm only when $1\leq p \leq \infty $. When $0<p<1$, $A^p(\Omega)$ becomes a metric space by defining $d(f,g)=\int_{\Omega} |f(z)-g(z)|^p \, dm(z)= \|f-g\|^p _{A^p(\Omega)}$ as usual.

\begin{theorem} \label{thm1.1}
The spaces $A^p(\Omega)$, $0<p<\infty$, are complete.
\end{theorem}

\begin{proof} It suffices to show that $A^p(\Omega)$ is a closed subspace of $L^p(\Omega)$. To this end, let $\{f_n\}$ be a sequence of functions in \ApO and $f\in L^p(\Omega)$ such that $\lim _{n\rightarrow \infty}\|f_n-f\|^p _{L^p(\Omega) }= 0$. Then there exists a subsequence $\{f_{n_k}\}$ of $\{f_n\}$ which converges to $f$ almost everywhere in $\Omega$. Moreover, since $\{f_n\}$ converges in $L^p(\Omega)$, $\{f_n\}$ is Cauchy in $L^p(\Omega)$. Due to Lemma \ref{lm1.1}, $\{f_n\}$ is uniformly Cauchy on the compact subsets of $\Omega$. Therefore $\{f_n\}$ converges uniformly on the compact subsets of $\Omega$ to some function $g$. Since each $f_n$ is analytic in $\Omega$, $g$ is also analytic in the same set. It turns out that $g$ must coincide with $f$ almost everywhere in $\Omega$ and hence that $f\in A^p(\Omega)$.    
\end{proof}

In what follows, we will use some classical notions from complex analysis as well as Montel's theorem. For the sake of completeness we shall digress a little and remind a few definitions as well as the statement of Montel's theorem.

\begin{definition}\label{def1.1} Let $G\subseteq \C$ be an open set and $(\Omega , d )$ a complete metric space. We define the space of continuous functions defined o $\Omega$ and taking their values in $G$ as
$$C(G,\Omega)=\{f:G\longrightarrow \Omega \; | \mbox{ f is continuous } \}.$$
\end{definition}

We usually consider the cases $\Omega=\C$ or $\Omega = \C_\infty $ in order to avoid trivial cases for the set $C(G,\Omega)$. For example, if $G$ is an open connected set of $ \C $ and $\Omega=\N=\{1,2,3,\ldots\}$, then $C(G,\Omega)$ contains just the constant functions on $G$.

The space $C(G,\C)$ given the metric of uniform convergence on the compact sets is a complete metric space. The space
$$H(G)=\{f:G\longrightarrow \C\; |\mbox{ f is analytic on G} \}\subseteq C(G,\C )$$
is also a complete metric space, if given the same metric.

\begin{definition}\label{def1.2} A set $\mathcal{F}\subseteq C(G,\Omega)$ is called \underline{normal} if every sequence of elements in $\mathcal{F}$ has a subsequence that converges in $C(G,\Omega)$ to some $f\in C(G,\Omega)$.
\end{definition}

We now state without proof Montel's theorem

\begin{theorem}\textnormal {(Montel)} \label{thm1.2}
Let $\mathcal{F}\subseteq H(G)$ be a family of functions. Then $\mathcal{F}$ is normal if and only if $\mathcal{F}$ is uniformly bounded on the compact subsets of $G$. 
\end{theorem}

Using Montel's theorem and the growth of $A^p$ functions we'll be able to show that every sequence in $A^p$, which is bounded with respect to the $A^p$ norm, has a subsequence that converges pointwise to some function $f\in A^p$.

\begin{proposition}\label{pr1.3} Let $\{f_n\}$ be a sequence of functions in $ A^p(\Omega)$ wich is uniformly bounded in $A^p(\Omega)$
 $$\|f_n\|_{A^p(\Omega)}\leq M \mbox{ for every }n\in\N.$$Then there exists a subsequence $\{f_{n_k}\}$ of $\{f_n\}$ and a function $f\in A^p(\Omega)$ such that $\{f_{n_k}\}$ converges to f uniformly on the compact subsets of $\Omega$. 
\end{proposition}
\begin{proof} From Proposition \ref{pr1.1} and the hypothesis we have that
$$|f_n(z)|\leq \pi^{-\frac{1}{p}} dist(z,\partial\Omega)^{-\frac{2}{p}}\|f_n\|_{A^p(\Omega)}\leq \pi^{-\frac{1}{p}} \delta^{-\frac{2}{p}}M$$
for every $z\in\Omega$. It is easy to see now that the sequence $\{f_n\}$ is uniformly bounded on the compact subsets of $\Omega$. According to Montel's theorem, this is equivalent to $\{f_n\}$ being normal on $\Omega$. By Definition \ref{def1.2}, this means that there exists a subsequence $\{f_{n_k}\}$ of $\{f_n\}$ and a function $f\in C(\Omega,\C)$ such that $\f_{n_k}$ converges to $f$ uniformly on the compact subsets of $\Omega$. It remains to show that $f$ is an element of $A^p(\Omega)$. But this is obvious since f must be analytic as a uniform limit of analytic functions.
\end{proof}

\noindent\textbf{Remark.} Baring in mind the growth condition (\ref{eq1.6})
\begin{equation*}
|f(z)|\leq \frac{1}{\pi ^{\frac{1}{p}}}\delta(z)^{-\frac{2}{p}}\|f\|_{A^p(\Omega)},
\end{equation*}
it is easy to see that $A^p(\C)=\{0\}$. In order to avoid trivial cases like this one we must be a little careful when chosing the set $\Omega$. In other words, $\Omega$ has to be an open connected set such that for every $z\in\Omega$ one can find an $f\in A^p(\Omega)$ satisfying $f(z)\neq 0$. There are however several open connected sets $\Omega$ which give rise to non trivial spaces $A^p(\Omega)$. We give some examples here.

\newpage

\noindent{(i)} The unit disc $\Delta$ and in general every bounded subset of $\C$.\\
\noindent{(ii)} Let $\Omega = \C -[0,\infty)$. Then the function $\psi:\Delta\longrightarrow \Omega$ defined as 
$$\psi(z)=\big(i\frac{1+z}{1-z} \big)^2 $$
 
is a conformal mapping (one to one, onto) satisfying $\psi(0)=-1$, $\psi(-1)=0$ and $\psi(1)=\infty$. Let $\zeta_o\in \Delta$. We want to find a function $f\in A^p(\Omega)$ satisfying $f(\zeta_o)\neq 0$. From (i), there exists a function $g\in A^p(\Delta)$ with $g(\psi^{-1}(\zeta_o))\neq0$. Define $f:\Omega \longrightarrow \C $ as $$f(w)=g(\psi^{-1})(w)=(g\circ \psi^{-1})(w).$$
Then 
\begin{eqnarray*}
\int_{\Omega}|f(\zeta)|^p dA(\zeta)&=&\int_{\Omega} |g\circ \psi^{-1}(\zeta)|^p dA(\zeta)=
\int_{\psi(\Delta)}|f(z)|^p|\psi^{'}(z)|^2 dA(z)\\
&=&\int_\Delta |f(z)|^p|\psi^{'}(z)|^2 dA(z)=4 \,\int_\Delta |f(z)|^p\bigg|i\frac{1+z}{1-z}\frac{2}{(1-z)^3} \bigg|^2 dA(z) \\
&=&16 \, \int_\Delta |f(z)|^p\bigg|\frac{1+z}{(1-z)^3} \bigg|^2 dA(z)<\infty.
\end{eqnarray*}
Since $f(\zeta_o)=g(\psi^{-1}(\zeta_o))\neq 0$, f is the function we were seeking for.

\pa In general, if $\Omega$ is an open, simply connected proper subset of $\C$ and $\zeta_o \in \Omega$, there always exists a unique one to one, analytic mapping $\phi$ from $\Omega$ onto $\Delta$ such that $\phi(\zeta _o)=0$ and $\phi^{'}(\zeta_o)>0$. Then
$$\int_{\Omega}|\phi^{'}(z)^\frac{2}{p}|^p dA(z)=\int_{\Omega}|\phi^{'}(z)|^2 dA(z)=\int_{\Delta} dA(z) <\infty.$$

\noindent{(iii)} Let $\Omega$ be an open connected. Suppose further that the  boundary of $\Omega$, $\partial \Omega$, has at least one connected component $S$ containing at least two points and that the the complement of $S$ in $\C_\infty$ is not a singleton. Then $\Omega$ defines a non trivial $A^p$ space.

\pa To see this observe first that since $\partial \Omega$ is closed and $S$ is a component, $S$ must also be closed and of course connected. Define $\Omega^*=\C_\infty-S\supseteq \Omega$. Clearly $\Omega^*$ is open and connected and therefore the discussion in (ii) yields that $A^p(\Omega^*)$ is non trivial. Since $\Omega\subseteq\Omega^*$ we are finished.

\medskip


\section{The Bergman Kernel}

\medskip

For this section, we'll restrict our attention to the Hilbert space case $p=2$. Let $\Omega$ be an open connected set such that the space \ApO is non trivial. Fix a $z\in \Omega$ and define the \textbf{point evaluation functional} 
$$\phi:A^2(\Omega)\longrightarrow \C, \hspace{5mm} \phi(f)=f(z).$$
The functional $\phi$ is bounded since for every $f\in A^2(\Omega)$ we have that
\begin{equation*}
|\phi(f)|=|f(z)|\leq \frac{1}{\pi ^{\frac{1}{p}}}\delta(z)^{-\frac{2}{p}}\|f\|_{A^2(\Omega)}
\end{equation*}
from Poroposition \ref{pr1.1} (remember that z is fixed). From the Riesz representation theorem for Hilbert spaces there exists a unique function $K_z\in A^2(\Omega)$ such that
$$\phi(f)=\langle f,K_z \rangle_{A^2(\Omega)}=\int_{\Omega}f(\zeta)\overline{K_z (\zeta)} dA(\zeta)$$
for every $f\in A^2(\Omega)$. 
\begin{definition}\label{def1.3} The function $K: \C\times \C \longrightarrow \C$ defined as 
$$K(z,\zeta)=\overline{K_z(\zeta)}$$
is called the \textbf{reproducing Kernel} or the \textbf{Bergman Kernel}.
\end{definition}
So the Bergman Kernel reproduces the values of every function in $A^2(\Omega)$ by means of the formula
\begin{equation} \label{eq1.9}
f(z)=\int_{\Omega} f(\zeta)K(z,\zeta)dA(\zeta).
\end{equation}

The following properties of the Bergman Kernel are simple consequences of the Definition \ref{def1.3}.

\noindent{(i)} For every $z,\zeta \in \C$ we have that
\begin{equation}\label{eq1.10}
K(z,\zeta)=\overline{K_z(\zeta)}=K_{\zeta}(z)=\overline{K(\zeta,z)}
\end{equation}

Indeed,
\begin{eqnarray*}
\overline{K(\zeta,z)}&=&K_{\zeta}(z)=\int_{\Omega} K_{\zeta}(w)K(z,w)dA(w)=
\int_{\Omega}\overline{K(\zeta,w)}K(z,w)dA(w)\\ \\
&=&\overline{\int_{\Omega}K(\zeta,w)\overline{K(z,w)}dA(w)}=
\overline{\int_{\Omega}K_z(w)	 K(\zeta,w)dA(w)}=\overline{K_z(\zeta)}\\
&=&K(z,\zeta).
\end{eqnarray*}

\noindent{(ii)} The function $K(z,\zeta)$ is analytic with respect to $z$ and counter analytic with respect to $\zeta$. This is just a consequence of (\ref{eq1.10}) since, by defintion, the functions $K_z$, $K_\zeta$ are in $A^2(\Omega)$ and hence analytic.

\noindent{(iii)} For every $z\in\Omega$ we have that 
\begin{equation}\label{eq1.11}
K(z,z)= \|K(z,\cdot ) \|^{2}_{A^2(\Omega)}.
\end{equation}

This is a just a simple calculation
\begin{eqnarray*}
K(z,z)&=& \overline{K(z,z)}=K_z(z)=\int_{\Omega}K_z(\zeta)\overline{K_z(\zeta)} dA(\zeta) 
\\ &=& \int_{\Omega}|K_z(\zeta)|^2= \|K(z,\cdot\|^{2}_{A^2(\Omega)}>0.
\end{eqnarray*}

\noindent{(iv)} Fix a $z\in\Omega$. Then, for every function $f\in A^2(\Omega)$ satisfying $f(z)=1$ we have the estimate
\begin{equation}\label{eq1.12} 
\|f\|_{A^2(\Omega)} \geq \|K(z,\cdot)\|^{-1}_{A^2(\Omega)} .
\end{equation}

To see this, use formula (\ref{eq1.9}) to write
\begin{eqnarray*}
|f(z)|&=&\bigg|\int_{\Omega}f(\zeta)K(z,\zeta)dA(\zeta)\bigg|\leq
\int_{\Omega}|f(\zeta)||K(z,\zeta)|dA(\zeta) \\ 
&\leq & \|f\|^2 _{A^2(\Omega)} \int_{\Omega} |K(z,\zeta)|^2dA(\zeta)=
 \|f\|^2 _{A^2(\Omega)}  \|K(z,\cdot)\|^2 _{A^2(\Omega)}.
\end{eqnarray*}
It is easy to see that estimate (\ref{eq1.12}) is sharp. Indeed, fix a $z\in\Omega$ and define 
$$f(\zeta)=\frac{K(\zeta,z)}{K(z,z)}=\frac{K(\zeta,z)}{\|K(z,\cdot)\|^2 _{A^2(\Omega)}}.$$
Then, $f(z)=1$ and 
\begin{eqnarray*}
\|f\|^2 _{A^2(\Omega)}&=&\frac{1}{\|K(z,\cdot)\|^2 _{A^2(\Omega)}} \int_{\Omega}|K(\zeta,z)|^2dA(\zeta) \\
 &=&\frac{\|K(z,\cdot)\|^2 _{A^2(\Omega)}}{\|K(z,\cdot)\|^4 _{A^2(\Omega)}}=\|K(z,\cdot)\|^{-2} _{A^2(\Omega)} .
\end{eqnarray*}

\noindent{(v)} The Bergman Kernel $K(z,\cdot)$ is the only function in $A^2(\Omega)$ which reproduces the value of every function $f\in A^2(\Omega)$ in the sense of (\ref{eq1.9}). This is just a consequence of the definition of the Bergman Kernel by means of the Riesz representation theorem.

\medskip
We next turn to the question of calculating the Bergman Kernel in $A^2(\Omega)$. We consider an orthonormal base $\{\phi_n\}^\infty _{n=0}$ of $A^2(\Omega)$. To simplify notation we write $\langle \cdot,\cdot \rangle$ for the inner product in $A^2(\Omega)$. Then
$$\langle\phi_n,\phi_m\rangle=\delta_{nm}=
\begin{cases}
                                             1,&\mbox{if} \hspace{4mm}n=m\;\\
                                            0   ,&\mbox{if} \hspace{4mm} n\neq m \; \ .\\      
                 \end{cases} $$

The usual Hilbert space identities hold:

\noindent(i) \textbf{"Fourier" series:} Every $f\in A^2(\Omega)$ has a series representation 
\begin{equation}\label{eq1.13}
f=\sum_{n=0} ^\infty c_n \phi_n, 
\end{equation}
where $c_n=\langle f,\phi_n \rangle$ and the series in (\ref{eq1.13}) converges with respect to the $A^2(\Omega)$ norm. Now, Lemma (\ref{lm1.1}) implies that the partial sums $\sum _{n=0} ^N c_n \phi_n$ converge to $f$ uniformly on the compact subsets of $\Omega$.

\noindent(ii) \textbf{Parseval's Identity:} For every in $f\in A^2(\Omega)$ we have that 
\begin{equation}\label{eq1.14}
\sum_{n=0} ^\infty |c_n|^2=\|f\|^2 _{A^2(\Omega)}.
\end{equation}

\noindent(iii) \textbf{A formula for the Bergman Kernel:} For every $z,\zeta \in \Omega$, the Bergman Kernel can be calculated as 
\begin{equation}\label{eq1.15}
K(z,\zeta)=\sum_{n=0} ^\infty \phi_n(z) \overline{\phi_n(\zeta)}
\end{equation}
where the series in (\ref{eq1.15}) converges with respect to the $A^2(\Omega)$ norm
and hence uniformly on the compact sets of $\Omega$.

To see this, write the series representation of the function $K(\cdot,\zeta)$ as in (\ref{eq1.13}). The "Fourier" coefficients of the function $K(\cdot,\zeta)$ are 
\begin{eqnarray*}
c_n&=&\langle K(\cdot,\zeta),\phi_n \rangle=\int_{\Omega} K(w,\zeta) \overline{\phi_n (w)}dA(w) =\overline{\int_{\Omega} \overline{K(w,\zeta)} \phi_n (w)dA(w)}\\
&=& \overline{\int_{\Omega} K(\zeta,w) \phi_n (w)dA(w)}=\overline{\phi_n(\zeta)}.
\end{eqnarray*}
Writing the series representation of the function $K(\cdot,\zeta)$
$$K(z,\zeta)=\sum^\infty _{n=0} c_n \phi_n(z)=\sum^\infty _{n=0}  \phi_n(z)\overline{\phi_n(\zeta)}$$
we get (\ref{eq1.15}).

\medskip
The next step is to try to calculate the Bergman Kernel in the case of the unit disc $\Omega=\Delta$. In order to do this we will employ formula (\ref{eq1.15}). First of all we have to define an orthonormal base in $A^2(\Delta)$.

\begin{proposition}\label{pr1.4} The set $ \{ \phi_n(z)=\sqrt{n+1} \; z^n\}_{n=0} ^\infty$ is an orthonormal base of $A^2(\Delta)$.
\end{proposition} 	
\begin{proof} The proof will be done in two steps:

\medskip
\noindent\textbf{step 1} The set $\{ \phi_n\}_{n=0} ^\infty$ is orthonormal.  

\smallskip
Indeed, for every $n,m\in\N$ we have that
\begin{eqnarray*}
\langle \phi_n , \phi _m \rangle &=&\sqrt{n+1}\sqrt{m+1}\int_{\Delta} z^n \overline{z^m} dm(z)
\\&=&
\sqrt{n+1}\sqrt{m+1}\int_0 ^1 \int_0 ^{2\pi}  (re^{i\theta})^n (re^{-i\theta})^m \frac{1}{\pi}d\theta rdr \\
&=&\frac{1}{\pi}\sqrt{n+1}\sqrt{m+1}\int_0 ^1 r^{n+m+1}dr \int _0 ^{2\pi} e^{i(n-m)\theta}d\theta
\\&=&\begin{cases}
    1\hspace{1mm}                        ,&\mbox{if} \hspace{4mm}n\ = m\;\\
    0\hspace{1mm}      ,&\mbox{if} \hspace{4mm} n\neq m \; \ . \\      
                 \end{cases} 
\end{eqnarray*}
\noindent\textbf{step 2} The set $\{ \phi_n\}_{n=0} ^\infty$ is a base of $A^2(\Delta)$.

\smallskip
It is equivalent to show that Parseval's formula
\begin{equation}\label{eq1.16}
\|f\|_{A^2(\Delta)}=\sum_{n=0} ^\infty |\langle f,\phi_n\rangle|^2
\end{equation}
holds true for every $f\in A^2(\Delta)$. If we write the expansion of $f$ in its Taylor series
$$f(z)=\sum_{n=0} ^\infty a_n z^n$$
it is clear that we have to show that
\begin{equation}\label{eq1.17}
\|f\|_{A^2(\Delta)}=\sum_{n=0} ^\infty \frac{|a_n|^2}{n+1} \; .
\end{equation}
Consider the partial sums of the Taylor series of $f$, $S_N(f)(z)=\sum_{n=0}^N a_nz^n$ and the disc $\Delta_\rho$ of radius $\rho$, centered at the origin, where $0<\rho<1$,
$$\Delta_\rho=\{z\in\C:|z|<\rho\}.$$ 
Then,
\begin{eqnarray*}\int_{\Delta_\rho}|S_N(f)(z)|^2dm(z)&=&\int_{\Delta_\rho}\bigg(\sum_{n=0}
^N a_nz^n\bigg)^2\bigg(\sum_{n=0}
^N \overline{a_n}\; \overline{z}^n\bigg)^2dm(z)\\
&=&
\sum_{n=0}^Na_n\sum_{m=0}^N \overline{a_n}  \int_0 ^\rho r^{n+m+1}\frac{1}{\pi}\int_0 ^{2\pi}e^{i(n-m)\theta}d\theta dr\\
&=&
2\sum_{n=0}^N|a_n|^2 \int_0 ^\rho r^{2n+1}dr=\sum_{n=0}^N|a_n|^2 \frac{\rho^{2(n+1)}}{n+1}\; .
\end{eqnarray*}
Since $S_N(f)$ converges to $f$ uniformly in $\Delta_\rho$ as $n\rightarrow \infty$, it is clear that
$$\lim_{N\rightarrow \infty} \int_{\Delta_\rho}|S_N(f)(z)|^2dm(z) =\int_{\Delta_\rho}|f(z)|^2dm(z).$$
However,
$$\lim_{N\rightarrow \infty} \sum_{n=0}^N|a_n|^2\frac{\rho^{2(n+1)}}{n+1}=\sum_{n=0}^\infty|a_n|^2\frac{\rho^{2(n+1)}}{n+1}$$
and so, combining the last two formulas we get
$$ \int_{\Delta_\rho}|f(z)|^2dm(z)= \sum_{n=0}^\infty|a_n|^2\frac{\rho^{2(n+1)}}{n+1}.$$
Letting $\rho\rightarrow 1^{-}$ yields (\ref{eq1.17}).
\end{proof}

Having defined an orthonormal base of $A^2(\Omega)$,it is now easy to calculate an exact formula for the Bergman Kernel. Indeed, employing (\ref{eq1.15}), we have that
$$K(z,\zeta)=\sum_{n=0} ^\infty \phi_n(z)\overline{\phi_n(\zeta)}=
\sum_{n=0} ^\infty(n+1)(z\overline{\zeta})^n =\frac{1}{(1-\overline{\zeta}z)^2} \; .$$

We have actually established the following.
\begin{theorem} \label{thm1.3}(i) The Bergman Kernel in $A^2(\Delta)$ is given by the formula
\begin{equation}\label{eq1.18}
K(z,\zeta)=\frac{1}{(1-\overline{\zeta}z)^2} \; .
\end{equation}
(ii) For every function $f\in A^2(\Delta)$ we have the representation
\begin{equation}\label{eq1.19}
f(z)=\int_{\Delta} \frac{f(\zeta)}{(1-\overline{\zeta}z)^2}\; dm(\zeta) \; .
\end{equation}
\end{theorem}

We close this section with a theorem that relates the Bergman Kernels of two open connected sets through a univalent mapping. 

\begin{theorem}\label{thm1.4}
Let $\Omega,\mathcal{D}$ be open and connected subsets of $\C$ and $\phi:\Omega\longrightarrow \mathcal{D}$, $\phi(z)=w$, be a univalent mapping of $\Omega$ onto $\mathcal{D}$. Suppose further that $J(w,\omega)$ is the Bergman Kernel in $\mathcal{D}$. Then, the Bergman Kernel in $\Omega$ is given by
\begin{equation}\label{eq1.20} 
K(z,\zeta)=J(\phi(z),\phi(\zeta))\phi^{'}(z)\overline{\phi^{'}(\zeta)}.
\end{equation}
\end{theorem}
\begin{proof}

Define the operator $T:A^2(\mathcal{D})\longrightarrow A^2(\Omega)$ as
$$T(f)(z)=(f\circ\phi)(z) \phi^{'}(z).$$
It easy to check that T is an isometry. Indeed,
\begin{eqnarray*}
\|f\|^2 _{A^2(\mathcal{D})}&=&\int_{\mathcal{D}}|f(w)|^2dA(w)=
\int_{\phi(\Omega)}|f(w)|^2dA(w)=\int_{\Omega}|f(\phi(z))|^2|\phi^{'}(z)|^2dA(w)\\
&=&\int_{\Omega}|T(f)(z)|^2dA(z)=\|Tf\|^2 _{A^2(\Omega)}.
\end{eqnarray*}
If $g\in A^2(\Omega)$ then define the function $f$ on $\mathcal{D}$ by the formula $$f(w)=g(\phi^{-1}(w))(\phi^{-1})^{'}(w).$$
Then $f\in A^2(\mathcal{D})$ by a simple change of variable and $T(f)=g$ which shows that $T$ is onto.

Now, consider any $g\in A^2(\Omega)$ and $f\in A^2(\mathcal{D})$ such thath $g=T(f)$. Then
$$f(w)=\int_{\mathcal{D}}J(w,\omega)f(\omega)dA(\omega).$$
and replacing $w$ by $\phi(z)$ in the above we get
$$f(\phi(z))=\int_{\mathcal{D}}J(\phi(z),\omega)f(\omega)dA(\omega)$$
for $z\in\Omega$. Making the change of variable $\omega=\phi(\zeta)$ results to
$$f(\phi(z))=\int_{\Omega}J(\phi(z),\phi(\zeta))f(\phi(\zeta))|\phi^{'}(z)|^2dA(\zeta)$$
and multiplying by $\phi^{'}(z)$,
$$f(\phi(z))\phi^{'}(z)=\int_{\Omega}J(\phi(z),\phi(\zeta))\phi^{'}(z)\overline{\phi^{'}(\zeta)}\ f(\phi(\zeta)) \phi^{'}(\zeta)dA(\zeta).$$
Remembering that $g(z)=T(f)(z)$ we get
$$g(z)=\int_{\Omega}J(\phi(z),\phi(\zeta))\phi^{'}(z)\overline{\phi^{'}(\zeta)}\ g(\zeta)dA(\zeta)$$
which is the desired result.
\end{proof}

An application of Theorem \ref{thm1.4} is contained in the Corollary below.

\begin{corollary}\label{cor1.1} Let $\Omega$ be an open, connected, proper subset of $\C$ and $K(z,\zeta)$ be the Bergman Kernel for $\Omega$. From Riemann's mapping theorem we know that for every $\zeta\in\Omega$ there exists a conformal mapping $\phi$ of $\Omega$ onto $\Delta$ with $\phi(\zeta)=0$ and $\phi^{'}(\zeta)>0$. Then, for every $z\in\Omega$, 
\begin{equation} \label{eq1.21}
\phi^{'}(z)=\sqrt{\frac{\pi}{K(\zeta,\zeta)}}K(z,\zeta) .
\end{equation}
\end{corollary}
\begin{proof}

Using theorem \ref{thm1.4}, the Bergman Kernel of $\Omega$ is written as
$$K(z,\zeta)=J(\phi(z),\phi(\zeta))\phi^{'}(z)\overline{\phi^{'}(\zeta)}$$
where $$J(w,\omega)=\frac{1}{\pi}\frac{1}{(1-w\overline{\omega})^2}$$ is the Bergman Kernel for the unit disc $\Delta$ (the constant $\frac{1}{\pi}$ is there because we have considered the normalised Lebesgue measure on the unit disc). Combining tha last two relations we get
$$K(z,\zeta)=\frac{1}{\pi}\frac{1}{(1-\phi(z)\overline{\phi(\zeta)})^2}\phi^{'}(z)\overline{\phi^{'}(\zeta)}=\frac{1}{\pi}\phi^{'}(z)\phi^{'}(\zeta)$$
since $\phi^{'}(\zeta)>0$ and $\phi(\zeta)=0.$ Therefore,
$$K(\zeta,\zeta)=\frac{1}{\pi}(\phi^{'}(\zeta))^2$$
and hence
$$K(z,\zeta)=\frac{1}{\pi}\phi^{'}(z)(\pi K(\zeta,\zeta))^\frac{1}{2}$$
which is just (\ref{eq1.21}). \end{proof}
\medskip


\section{Some Density Matters}

\medskip

It is obvious from Proposition \ref{pr1.4} that the polynomials are dense in $A^2(\Delta)$. In fact, something stronger is true.

\begin{theorem}\label{thm1.5} The set of polynomials is dense in $A^2(\Delta)$. Whatsmore, every function in $A^2(\Delta)$ can be approached in the $A^2(\Delta)$ norm by the partial sums of its Taylor series, that is
\begin{equation}\label{eq1.22}
 \lim_{N\rightarrow \infty}\|S_N(f)-f\|_{A^2(\Delta)}=0
\end{equation}
where $S_N(f)=\sum_{n=0} ^N a_n z^n$ is the partial sum of the Taylor series of $f$.
\end{theorem}

In order to prove this we shall need a simple lemma that relates the Taylor coefficients of $f$ with its "Fourier" coefficients with respect to the orthonormal base $\{\phi_n\}$.

\begin{lemma}\label{lm1.2} Suppose that $f\in A^2(\Delta)$ has the Taylor expansion $f(z)=\sum_{n=0}^\infty a_nz^n$ and that $\phi_n= \sqrt{n+1}\; z^n$. Then, for every $n\in \N$, 
\begin{equation}\label{eq1.23}
\langle f,\phi_n \rangle =\frac{a_n}{\sqrt{n+1}}
\end{equation}
\end{lemma}
\begin{proof} Consider the disc $\Delta_\rho=\{z\in\C:|z|<\rho\}$ and fix some $n\in \N$. Then, for $N>n$ we have that
\begin{eqnarray*}
\int_{\Delta_\rho}S_N(f)(z)\overline{\phi_n(z)}dm(z)&=&\int_{\Delta_\rho}\bigg(\sum_{k=0} ^Na_kz^k(z)\bigg)\overline{\phi_n(z)}dm(z) \\
&=&\sqrt{n+1}\sum_{k=0} ^Na_k  \int_{\Delta_\rho} z^k\overline{z}^n dm(z) \\
&=&\sqrt{n+1}\ 2 \ a_n\int_0 ^\rho r^{2n+1}dr=\sqrt{n+1}a_n \frac{\rho^{2n+1}}{n+1} \; .
\end{eqnarray*}
Since $S_N(f)$ converges to $f$, uniformly on $\Delta_\rho$ as $N\rightarrow \infty$, we get
$$\int_{\Delta_\rho}S_N(f)(z)\overline{\phi_n(z)}dm(z)=\sqrt{n+1}a_n \frac{\rho^{2n+1}}{n+1}\ .$$
Letting $\rho\rightarrow 1^{-}$ we get (\ref{eq1.23}).
\end{proof}
\noindent\textit{Proof of Theorem \ref{thm1.5}.} We write the Taylor series of f as
$$f(z)=\sum_{n=0} ^\infty a_n z^n=\sum_{n=0} ^\infty \frac{a_n}{\sqrt{n+1}}\sqrt{n+1}z^n=\sum_{n=0} ^\infty \frac{a_n}{\sqrt{n+1}} \ \phi_n(z)\ .$$
Then, for $N\in\N $,
\begin{eqnarray*}
&&\int_{\Delta}|f(z)-S_N(f)(z)|^2dm(z)=\langle f-S_N(f),f-S_N(f) \rangle
\\&=&\|f\|^2 _{A^2(\Delta)}-\langle f,S_N(f)\rangle-\overline{\langle f,S_N(f)\rangle}+\|S_N(f)\|^2 _{A^2(\Delta)}\\
&=&\|f\|^2 _{A^2(\Delta)}-\sum_{n=0} ^N \overline{a_n}\langle f,z^n\rangle-\sum_{n=0}^N a_n\overline{\langle f,z^n\rangle}+\sum_{n=0} ^N \frac{|a_n|^2}{n+1} \\ 
&=&\|f\|^2 _{A^2(\Delta)}-\sum_{n=0} ^N \frac{\overline{a_n}}{\sqrt{n+1}}\langle f,\phi_n\rangle-\sum_{n=0}^N \frac{a_n}{\sqrt{n+1}}\overline{\langle f,\phi_n\rangle}+\sum_{n=0} ^N \frac{|a_n|^2}{n+1} \ .
\end{eqnarray*}	
Employing Lemma \ref{lm1.2}, the right hand side is equal to
\begin{eqnarray*}
\|f\|^2 _{A^2(\Delta)}-\sum_{n=0} ^N \frac{|a_n|^2}{n+1}-\sum_{n=0}^N \frac{|a_n|^2}{n+1}+\sum_{n=0} ^N \frac{|a_n|^2}{n+1}=
\|f\|^2 _{A^2(\Delta)}-\sum_{n=0} ^N \frac{|a_n|^2}{n+1} \ .
\end{eqnarray*}
However, $$\|f\|^2 _{A^2(\Delta)}=\sum_{n=0} ^\infty \frac{|a_n|^2}{n+1}$$
from Parseval's identity and so
$$\|f-S_N(f)\|_{A^2(\Omega)} \longrightarrow 0\hspace{3mm} \mbox{  as  }\hspace{3mm} N\rightarrow \infty$$
and this finishes the proof.\qed

\smallskip
\noindent\textbf{Remark.} It is \underline{not} true that the polynomials are dense in $A^2(\Omega)$ for every simply connected set $\Omega\subseteq\C$. Consider for example the set $\Omega=\Delta-[0,1]$. This is obviously simply connected and gives rise to the space $A^2(\Omega)$. On $\Omega$ one can consider the function $f(z)=z^\frac{1}{2}$ which is clearly analytic square integrable. However, it is not possible to approach $f$ by polynomials in the $A^2(\Omega)$ sense.

To see this, suppose that one could find a sequence of polynomials $\{P_n\}_{n\in\N}$, such that
$$\lim_{n\rightarrow \infty}\|f-P_n\|_{A^2(\Omega)}=0\ . $$
Then the sequence $\{P_n\}_{n\in\N}$ is Cauchy in $A^2(\Omega)$ and hence in $A^2(\Delta)$ since $\|P_n-P_m\|_{A^2(\Omega)}=\|P_n-P_m\|_{A^2(\Delta)}$. This means that the sequence $\{P_n\}_{n\in\N}$ is uniformly Cauchy on the compact subsets of $\Delta$ which in turn means that it converges uniformly on the compact subsets of $\Delta$ to some $A^2(\Delta)$ function, say $g$. Since $\{P_n\}_{n\in\N}$ also converges to $f$ uniformly on the compact subsets of $\Omega$, it turns out that $f\equiv g$ in $\Omega$. But this means that $f$ has an analytic extension to the whole of $\Delta$, a contradiction.
\medskip

We will be able to show next that the polynomials are dense in $A^p(\Delta)$ for general $p$, $0<p<\infty$. However, we wont be able to approach a function in $A^p(\Delta)$ by the partial sums of its Taylor series.

\begin{theorem}\label{thm1.6}
The polynomials are dense in $A^p(\Delta)$, $0<p<\infty$.
\end{theorem}
\begin{proof} Let $f\in A^p(\Delta)$. We consider the function $f_\rho=f(\rho z)$, $0<\rho<1$. The function $f_\rho$ is analytic inside the disc $\Delta_\rho=\{z\in\C:|z|<\frac{1}{\rho}\}\supset \overline{\Delta}$. We deduce that the partial sums of the function $f_\rho$, $S_N(f_\rho)$, converge to $f_\rho$ uniformly on the compact subsets of $\Delta_\rho$ and hence uniformly in $\Delta$. That is
\begin{equation}\label{eq1.24}
S_N(f_\rho)_{}\longrightarrow f_\rho \mbox{\hspace{2mm} uniformly in\hspace{2mm} } \Delta \hspace{2mm}\mbox{ as }\hspace{2mm} N\rightarrow \infty.
\end{equation}  
It is immidiate from (\ref{eq1.24}) that
\begin{equation}\label{eq1.25}
\lim_{N\rightarrow \infty}\|S_N(f_\rho)-f_\rho\|^p _{A^p(\Delta)}=0.
\end{equation} 
It suffices to show that 
\begin{equation}\label{eq1.26}
\lim_{\rho\rightarrow 1}\|f-f_\rho\|^p _{A^p(\Delta)}=0.
\end{equation} 
Indeed, assuming for a moment (\ref{eq1.25}), we have that
$$\|f-S_N(f_\rho)\|^p _{A^p(\Delta)}\leq 2^p\|f-f_\rho\|_{A^p(\Delta)}+ 2^p\|f_\rho-S_N(f_\rho))\|_{A^p(\Delta)}.$$

Let $\epsilon>0$. We chose $\rho$ close to 1 so that $2^p\|f-f_\rho\|_{A^p(\Delta)} < \frac{\epsilon}{2}$ (this is possible because of (\ref{eq1.26})). Then, for $N$ large enough, $2^p\|f_\rho-S_N(f_\rho))\|_{A^p(\Delta)}<\frac{\epsilon}{2}$ because of (\ref{eq1.25}) and so
$$\|f-S_N(f_\rho)\|^p _{A^p(\Delta)} < \epsilon $$
and so $S_N(f_\rho)$ is the seeked for polynomial.

It remains to prove equation (\ref{eq1.26}). We have that

\begin{eqnarray*}
\|f-f_\rho\|^p _{A^p(\Delta)}&=&\int_\Delta|f(z)-f_\rho(z)|^p dm(z)
\\&=&\int_0 ^1 \frac{1}{\pi} \int _0 ^{2\pi} |f(re^{i\theta})-f_\rho(re^{i\theta})|^pd\theta r dr
\end{eqnarray*}

For the inner integral we have the estimate

\begin{eqnarray*}
\int _0 ^{2\pi} |f(re^{i\theta})-f_\rho(re^{i\theta})|^pd\theta &\leq&
\int _0 ^{2\pi} (|f(re^{i\theta})|+|f_\rho(re^{i\theta})|)d\theta  \\
&\leq & 2^p\int _0 ^{2\pi} \bigg\{|f(re^{i\theta})|^p+|f_\rho(re^{i\theta})|^p \bigg\}d\theta \\
&=& 2^p \bigg\{ \int _0 ^{2\pi} |f(re^{i\theta})|^pd\theta +  \int _0 ^{2\pi} |f(\rho re^{i\theta})|^p d\theta \bigg\}\\
&\leq&2^{p+1}  \int _0 ^{2\pi} |f(re^{i\theta})|^pd\theta<\infty .
\end{eqnarray*}

The last inequality is due to the fact the function
$$F(r)=\int_0 ^{2\pi}|f(re^{i\theta})|^p d\theta, \;\;\;0\leq r<1,$$
is an increasing function of r. Thus, since $0<\rho r<1$ ($0<\rho<1$), we have that $$F(\rho r )<F(r).$$ This is a simple lemma which we'll prove after the proof of this theorem. Thus we get
\begin{equation}\label{eq1.27}
\int _0 ^{2\pi} |f(re^{i\theta})-f_\rho(re^{i\theta})|^pd\theta \leq 2^{p+1}  \int _0 ^{2\pi} |f(re^{i\theta})|^pd\theta
\end{equation}

We show next that $f_\rho$ converges to $f$ uniformly on the compact subsets of $\Delta$ as $\rho \rightarrow 1$. If $K$ is a compact subset of $\Delta$ and $z\in K$ then 
$$|f(z)-f\rho(z)|=|\sum_{n=0} ^\infty a_n(1-\rho^n)z^n|\leq \sum_{n=0} ^\infty |a_n||1-\rho^n||z|^n\leq \sum_{n=0} ^\infty |a_n||1-\rho^n|M^n$$
for some $M<1$. Now, Lebesgue's dominated convergence theorem for series yields that
$f_\rho$ converges to $f$ uniformly on the compact subsets of $\Delta$ as $\rho \rightarrow 1$. Thus, for $r<1$ fixed, we get
$$\lim_{\rho \rightarrow 1} |f(re^{i\theta})-f_\rho (re^{i\theta})|=0$$
and so
\begin{equation}\label{eq1.28}
\lim_{\rho \rightarrow 1} \int_0 ^{2\pi}|f(re^{i\theta})-f_\rho (re^{i\theta})|d\theta=0.
\end{equation}
Combining (\ref{eq1.27}) with (\ref{eq1.28}) and Lebesgue's dominated convergence theorem we get (\ref{eq1.26}).
\end{proof}

We now give the proof of the Lemma we already used in the proof of Theorem \ref{thm1.6}.

\begin{lemma}\label{lm1.3} Suppose $g$ is a nonnegative subharmonic function in $\Delta$ and $0\leq r <1$. Then the function
$$F(r)=\frac{1}{\pi}\int_0 ^{2\pi}g(re^{i\theta}) d\theta.$$
is a decreasing function of $r$.
\end{lemma}
\begin{proof} Since $g$ is a subharmonic function, for every open connected set $B$ with $\overline{B}\subset \Delta$, there exists a function $U$, harmonic in $B$, such that
$g(z)= U(z)$ in $\partial B$ and $g(z)\leq U(z)$ in $B$.

Suppose that $0\leq r_1 <r_2<1.$ Set $B=\{z\in \C:|z|<r_2\}$. Then

$$F(r_1)=\frac{1}{\pi}\int_0 ^{2\pi} g(r_1e^{i\theta}) d\theta \leq 
\frac{1}{\pi}\int_0 ^{2\pi} U(r_1 e^{i\theta}) d\theta =2U(0)$$
by the mean value theorem. Again by the mean value theorem we have that $$2U(0)=\frac{1}{\pi}\int_0 ^{2\pi}U(r_2 e^{i\theta})d\theta $$
and so
$$F(r_1)\leq\frac{1}{\pi}\int_0 ^{2\pi}U(r_2 e^{i\theta})d\theta = \frac{1}{\pi}\int_0 ^{2\pi}
g(r_2e^{i\theta}) d\theta =F(r_2)$$
which shows that $F$ is increasing.
\end{proof}

\noindent\textbf{Remark.} In theorem \ref{thm1.6} we used the fact that the function
$$F(r)=\int_0 ^{2\pi}|f(re^{i\theta})|^p d\theta, \;\;\;0\leq r<1,$$
is an increasing function of $r$. This is an immidiate consequence of lemma \ref{lm1.3} since the function $|f(z)|^p$ is a subharmonic function.


\chapter {The Bergman Projection }
\label{ch2}
 
\section {The Bergman Projection on $A^p(\Delta)$}
\medskip

Let us recall the formula of the main theorem of the previous chapter, that is 
Theorem \ref{thm1.3}. The latter states that for every function $f\in A^2(\Delta)$ we have the representation
\begin{equation}\label{eq2.1} 
f(z)=\int_\Delta \frac{f(\zeta)}{(1-\overline{\zeta}z)^2}\ dm(\zeta).
\end{equation}
Although the discussion that led to Theorem \ref{thm1.3} strongly depends on the Hilbert space structre of $A^2(\Delta)$, one can try to see if the above formula has a meaning in a more general context. This is the content of the following proposition.

\begin{proposition}\label{pr2.1}
For $f\in L^p(\Delta)$, $1\leq p < \infty $, define the function $P(f)$ on $\Delta$ as
\begin{equation}\label{eq2.2}
P(f)(z)=\int_\Delta \frac{f(\zeta)}{(1-\overline{\zeta}z)^2}\ dm(\zeta), \ \ \ z\in \Delta .
\end{equation}
Then,\\
\noindent(i) The function $P(f)$ is a well defined analytic function on $\Delta$.\\
\noindent(ii) If in addition $f\in A^p(\Delta)$, $1\leq p <\infty$ and $z\in \Delta$, we have that
\begin{equation}\label{eq2.3}
P(f)(z)=f(z).
\end{equation}
\end{proposition} 
\begin{proof}
Let us first note that the integral in (\ref{eq2.2}) is well defined for every $f\in L^p(\Delta)$, $1\leq p < \infty$. Indeed, fix a $z\in \Delta $ and let $1\leq p <\infty$ and $q$ be the conjugate exponent of $p$, that is, $\frac{1}{p}+\frac{1}{q}=1$. Then,
$$\int_\Delta \frac{|f(\zeta)|}{|1-\overline{\zeta}z|^2}\ dm(\zeta) \ \leq 
\Big\{\int_{\Delta} \frac{1}{|1-\overline{\zeta}z|^{2q}} \ dm(\zeta)\Big\}^\frac{1}{q}	\|f\|_{L^p(\Delta)}.$$
However, for $z$ fixed and $\zeta \in \Delta $, the function $\frac{1}{(1-\overline{\zeta}z)^2}$ is a bounded function of $\zeta $ and hence in every $L^q(\Delta)$, $1<q\leq \infty$. Since the integral in (\ref{eq2.2}) defines an analytic function on $\Delta$ whenever it exists, this proves (i).

For (ii) observe that the spaces $A^p(\Delta)$ are nested so Theorem \ref{thm1.3} implies that formula (\ref{eq2.3}) holds true for every function $f\in A^p(\Delta)$, $2 \leq p<\infty $.
It's easy to extend this formula to $A^1(\Delta)$ and hence to every $A^p(\Delta)$, $1\leq p <\infty $. Indeed, consider an $f\in A^1(\Delta)$. For $\rho \in (0,1)$ define the function $f_\rho(z)=f(\rho z)$. Now, $f_\rho \in H^\infty(\Delta)\subset A^2(\Delta)$ and therefore we can write
$$f_\rho(z)=\int_\Delta \frac{f_\rho(\zeta)}{(1-\overline{\zeta}z)^2}dm(\zeta)
=\frac{1}{\rho}\int_{\Delta} \frac{f(\zeta)}{(1-\frac{\overline{\zeta}}{\rho} z)^2} \chi _{\Delta_\rho}(\zeta)\ d\zeta $$ 
where $\chi_{\Delta_\rho}$ is the characteristic function of the disc of radius $\rho$ centered at the origin. For every $\zeta \in \Delta$ and $\rho \in (0,1)$ we have that
$$\frac{|f(\zeta)|}{|1-\frac{ \overline{\zeta}}{\rho}z|^2} \chi _{\Delta_\rho}(\zeta) 
\leq \frac{1}{|1-z|^2}|f(\zeta)|\in L^1(\Delta).$$
Empolying Lebesgue's dominated convergence theorem yields formula (\ref{eq2.3}) for $f\in A^1(\Delta)$. 
\end{proof}
\begin{definition}\label{def2.1} The linear operator P  is called the \textbf{Bergman Projection}.
\end{definition}

In the stronger $L^2(\Delta)$ case, one can easily see that the Bergman Projection is the orthogonal projection of $L^2(\Delta)$ onto the closed subspace $A^2(\Delta)$.  

\begin{proposition}\label{pr2.2} The Bergman Projection is the orthogonal projection of $ L^2(\Delta)$ onto $A^2(\Delta)$.
\end{proposition}
\begin{proof} Since $A^2(\Delta)$ is a closed subspace of $L^2(\Delta)$, there exists an orthogonal projection, say $P_o$, of $L^2(\Delta)$ onto $A^2(\Delta)$. Then, for every $f\in L^2(\Delta)$, $P_o(f)\in A^2(\Delta)$ and so
$$P_o(f)(z)=\langle P_o(f),K_z \rangle = \langle f,P_o(K_z) \rangle =\langle f,K_z \rangle =P(f).$$ 
Consequently, $P_o$ coincides with $P$ which shows the proposition.
\end{proof}

When $p\neq 2$ there is no orthogonal projection. However, since $P(L^p(\Delta))\supseteq A^p(\Delta)$, it is natural to ask whether
the Bergman Projection is a bounded operator from $L^p(\Delta)$ to $A^p(\Delta)$ which would also show that $P(L^p(\Delta))= A^p(\Delta)$.

Let us first show that this is not the case when $p=1$. 

\begin{proposition}\label{pr2.3} The Bergman projection is not bounded from $L^1(\Delta)$ to $L^1(\Delta)$.
\end{proposition}
\begin{proof} We will actually show that the adjoint operator of $P$, $P^*$, is not a bounded operator on $L^\infty (\Delta)$. We therefore need to find a formula for the adjoint operator. To that end, consider $f\in L^1(\Delta)$ and $h\in L^\infty(\Delta)$. Writing down the definition of $P^*$ we have that
\begin{eqnarray*}
\langle f,P^*(h)\rangle&=&\langle P(f),h\rangle=\int_{\Delta}P(f)(z)\overline{h(z)}dm(z)=
\int_\Delta\Bigg\{ \int_\Delta \frac{f(\zeta)}{(1-\overline{\zeta}z)^2}dm(\zeta)\Bigg\} \overline{h(z)}\ dm(z) \\
&=&\int_\Delta f(\zeta)\overline{ \Bigg\{ \int_\Delta \frac{h(z)}{(1-\overline{z}\zeta)^2}dm(z)\Bigg\}}\ dm(\zeta),
\end{eqnarray*}
where the last equality follows by applying Fubini's theorem. Thus, we have established the formula
$$P^*(h)=\int_\Delta \frac{h(\zeta)}{(1-\overline{\zeta} z)^2}\ dm(\zeta), \ \ \ h\in L^\infty(\Delta).$$

Suppose now, for the sake of contradiction, that $P$ is a bounded operator on $L^1(\Delta)$, that is, that $P\in B(L^1(\Delta))$\footnote{We denote by $B(L)$ the set of bounded linear operators from $L$ to $L$.}. This is equivalent to saying that $P^*\in B(L^\infty(\Delta))$. For $z\in \Delta$ and $a\in (0,1)$ define the functions $g_a$ as 
$$g_a(\zeta)=\frac{(1-a\overline{\zeta})^2}{|1-a\overline{\zeta}|^2}.$$
Clearly $g_a \in L^\infty(\Delta)$ and $\|g_a\|_{L^\infty(\Delta)}=1$. However, 
\begin{eqnarray*}
P(g_a)(a)&=&\int_\Delta \frac{g_a(\zeta)}{(1-\overline{\zeta} z)^2}\ dm(\zeta)=
\int_\Delta \frac{1}{|1-a\overline{\zeta}|^2}\ dm(\zeta)
\\&=&2\int_0 ^1 \frac{1}{2\pi}\int _0 ^{2\pi} \frac{1}{|1-are^{-i\theta}|^2}\ d\theta \ r dr =2\int_0 ^1 \sum _{n=0} ^\infty a^{2n} r^{2n} r dr \\
&=&2\int_0 ^1  \sum _{n=0} ^\infty a^{2n} r^{2n+1} dr = \sum_{n=0} ^\infty \frac{1}{n+1}a^{2n}=\log{\frac{1}{1-a^2}}.
\end{eqnarray*}
Now the hypothesis that $P^* \in B(L^\infty)$ implies that
$$\log{\frac{1}{1-|a|^2}}\leq \|P^* g_a\|_{L^\infty(\Delta)}\leq \|g_a\|_{L^\infty(\Delta)}=1 .$$
Since this must hold for every $a\in(0,1)$ we ger a contradiction as $a\rightarrow 1^-$.
\end{proof}

Having got rided of the "bad" case $p=1$ we will now show that the Bergman projection is indeed a bounded operator from $L^p(\Delta)$ onto $A^p(\Delta)$ for all $1<p<\infty$. 

\begin{theorem}\label{thm2.1} The Bergman projection is a bounded linear operator from $L^p(\Delta)$ onto $A^p(\Delta)$ for every $1<p<\infty$.
\end{theorem}
\begin{proof} We have already showed that $P$ is onto since $P(f)=f$  for every $f\in A^p(\Delta)\subset L^p(\Delta)$. Now, fix $1<p<\infty$ and an $f\in L^p(\Delta)$. It is clear that $P(f)$ defines an anlytic function in the unit disc so it remains to show that $P(f)$ is in $L^p(\Delta)$. For $z\in \Delta$ and $q$ the conjugate exponent of $p$ we have the estimate
\begin{eqnarray*}
|P(f)(z)|&\leq& \int_\Delta \frac{|f(\zeta)|}{|1-\overline{\zeta}z|^2}\ dm(\zeta)
=\int_\Delta \frac{(1-|\zeta|^2)^{-\frac{1}{pq}}}{|1-\overline{\zeta}z|^\frac{2}{q}}          \frac{|f(\zeta)|(1-|\zeta|^2)^\frac{1}{pq}}{|1-\overline{\zeta}z|^\frac{2}{p}}\ dm(\zeta)
\\ &\leq& \Bigg\{\int_\Delta \frac{(1-|\zeta|^2)^{-\frac{1}{p}}}{|1-\overline{\zeta}z|^2}\ \ dm(\zeta) \Bigg\}^\frac{1}{q}
  \Bigg\{\int_\Delta   \frac{(1-|\zeta|^2)^\frac{1}{q}}{|1-\overline{\zeta}z|^2} |f(\zeta)|^p \ dm(\zeta) \Bigg\}^\frac{1}{p} \\
  &=&[ I_1(z)]^\frac{1}{q}[ I_2(z)]^\frac{1}{p}.
\end{eqnarray*}
Taking $p$\textsuperscript{\textit{th}} powers and integrating we get
\begin{equation}\label{eq2.4}
\int_\Delta |P(f)(z)|^p \ dm(z) \leq \int_\Delta [I_1(z)]^\frac{p}{q} I_2(z) \ dm(z).
\end{equation}

At this point we need to estimate integrals of the form $\int_\Delta \frac{(1-|\zeta|^2)^\alpha}{|1-\overline{\zeta}z|^\beta} \ dm(\zeta)$ for suitable choices of $\alpha,\beta \in \R$. Instead of doing so in the special case we're interested in, we will state and prove a general lemma which will come in handy in many cases through-out the text.

\begin{lemma}\label{lm2.1} Let $s,t \in \R$ with $1<t<s$. Then there exists a constant $C=C(s,t)$, depending only on $s,t$, such that
\begin{equation}\label{eq2.5}
\int_\Delta \frac{(1-|\zeta|^2)^{t-2}}{|1-\overline{z}\zeta|^s} \ dm(\zeta) \leq C(1-|z|^2)^{t-s},
\end{equation}
for every $z\in\Delta$.
\end{lemma}

Postponing the proof of this lemma for a while, let's see how we can apply it to complete the proof of the theorem. For $I_1$ take $t=2-\frac{1}{p}$ and $s=2$ in the lemma above to get $I_1(z)\leq C(1-|z|^2)^{-\frac{1}{p}}.$ Plugging this estimate into formula (\ref{eq2.4}) we get
\begin{eqnarray*}
\int_\Delta |P(f)(z)|^p \ dm(z) &\leq & C \int_\Delta (1-|z|^2)^{-\frac{1}{q}}\Bigg\{\int_\Delta 
\frac{(1-|\zeta|^2)^\frac{1}{q}}{|1-\overline{\zeta}z|^2}|f(\zeta)|^p\ dm(\zeta) \Bigg\}
\ dm(z) \\
&=& C \int_\Delta |f(\zeta)|^p(1-|\zeta|^2)^\frac{1}{q} \Bigg\{\int _\Delta \frac {(1-|z|^2)^{-\frac{1}{q}}}{|1-\overline{\zeta}z|^2} \ dm(z) \Bigg\} \ dm(\zeta),
\end{eqnarray*}
where the last equation comes from an application of Fubini's theorem. Now, using Lemma \ref{lm2.1} on more time with $t=2-\frac{1}{q}$ and $s=2$ we  get
$$\int_\Delta |P(f)(z)|^p \ dm(z)\leq C^\prime  \int_\Delta |f(\zeta)|^p (1-|\zeta|^2)^\frac{1}{q}(1-|\zeta|^2)^{-\frac{1}{q}}\ dm(\zeta)=C^\prime \|f\|^p _{L^(\Delta)}$$
which finishes the proof.
\end{proof}

We now prove Lemma \ref{lm2.1}.

\medskip
\noindent\textit{Proof of Lemma \ref{lm2.1}.} Let $s,t\in \R$ with $1<t<s$ and set $I=\int_\Delta \frac{(1-|\zeta|^2)^{t-2}}{|1-\overline{z}\zeta|^s} \ dm(\zeta)$. Using the fact that $\frac{1}{(1-\overline{z}\zeta)^\frac{s}{2}}=\sum_{n=0}^\infty \frac{\Gamma (n+\frac{s}{2})}{n!\Gamma(\frac{s}{2})}\ \overline{z}^n \zeta^n$ and writing $I$ in polar coordinates we have
\begin{eqnarray*} 
I&=&\int_0 ^1 \frac{1}{\pi} \int _0 ^{2\pi} \frac{1}{|(1-\overline{z}re^{i\theta})^\frac{s}{2}|^2}d\theta(1-r^2)^{t-2}r \ dr \\
&=&2\sum _{n=0} ^\infty \Bigg(\frac{\Gamma (n+\frac{s}{2})}{n!\Gamma(\frac{s}{2})} \Bigg)^2\int_0 ^1 r^{2n+1}(1-r^2)^{t-2} \ dr |z|^{2n} \\
&=&\sum _{n=0} ^\infty \Bigg(\frac{\Gamma (n+\frac{s}{2})}{n!\Gamma(\frac{s}{2})} \Bigg)^2\int_0 ^1 r^n(1-r)^{t-2} \ dr |z|^{2n}\\
&=&\sum _{n=0} ^\infty \Bigg(\frac{\Gamma (n+\frac{s}{2})}{n!\Gamma(\frac{s}{2})} \Bigg)^2\frac{\Gamma(n+1)\Gamma(t-1)}{\Gamma(n+t)}|z|^{2n}\\
&=&\frac{\Gamma(t-1)}{\Gamma(\frac{s}{2})^2}\sum _{n=0} ^\infty \Bigg(\frac{\Gamma (n+\frac{s}{2})}{n!} \Bigg)^2\frac{n!}{\Gamma(n+t)}|z|^{2n}.
\end{eqnarray*}
Using standard estimates for the Gamma function we see that $ \Big(\frac{\Gamma (n+\frac{s}{2})}{n!} \Big)^2\frac{n!}{\Gamma(n+t)}$	is of the order $\frac{\Gamma(n+s-t)}{n!}$. But this means that
$$I\leq C(s,t) \sum_{n=0} ^\infty \frac{\Gamma(n+s-t)}{n!}|z|^{2n} =C(1-|z|^2)^{t-s}$$
which is just the statement of the lemma.\qed

\medskip


\section {A Bounded Projection of $L^1(\Delta)$ onto $A^1(\Delta)$}
\medskip

The Bergman projection seems to be the "natural" operator that maps $L^p(\Delta)$ onto  $A^p(\Delta)$, that is, it defines for every $L^p(\Delta)$ function, its "analytic counterpart" on the unit disc and it is the identity when restricted to $A^p(\Delta)$. The image of $L^p(\Delta)$ under the Bergman projection is exactly $A^p(\Delta)$ which reflects the fact that, well, it's a projection! This nice theory however, fails to provide with a bounded projection of $L^1(\Delta)$ onto $A^1(\Delta)$ since all these nice properties hold for $1<p<\infty$. We therefore set ourselves the task to find out if such a projection exists on $L^1(\Delta)$ and, if it does, to describe it. As a sideresult, we will also define the weighted Bergman spaces which arise naturally in the seek of such a projection.

First, recall the representation formula
$$f(z)=\int_\Delta\frac{f(\zeta)}{(1-\overline{\zeta}z)^2} \ dm(\zeta),\  \ z\in \Delta,$$
for a "nice" function $f$, say $f\in H^\infty(\Delta)$. Although this identiy also holds for $A^1$ functions, it does not define a bounded operator. As a first step, we will try to construct a family of representation formulas, at least for nice functions. This family will be more general, in the sense that it will include the above as a special case. Then we will see under what hypothesis this new family may define a bounded projection of $L^1(\Delta)$ onto $A^1(\Delta)$.

A first remark is that it suffices to represent a specific value of a function $f$, say $f(0)$. It is then easy to use a disc automorphism that carries 0 to any $z\in \Delta$ and automatically obtain a representation formula for the values of $f$ at any $z\in\Delta$. In this spirit, we have the following lemmas.

\begin{lemma}\label{lm2.2}
Let $f\in H^\infty(\Delta)$ and $\alpha>-1$. Then
\begin{equation}\label{eq2.6}
f(0)=(\alpha+1)\int_\Delta f(\zeta)(1-|\zeta|^2)^\alpha \ dm(\zeta).
\end{equation}
\end{lemma}
\begin{proof} 
First of all notice that the integral is well defined since $f$ is bounded and $\alpha>-1$. Now, if $f(z)=\sum_{n=0} ^\infty a_n z^n $ is the Taylor series of $f$, the right hand side of (\ref{eq2.6}) equals
\begin{eqnarray*}
	  \int  _\Delta  \sum_{n=0} ^\infty a_n \zeta^n (1-|\zeta|^2)^\alpha \ dm(\zeta)
&=&\int_0 ^1 \bigg\{\sum_{n=0} ^\infty \frac{1}{\pi} \int_0 ^{2\pi} e^{in\theta} d\theta \ r^n\bigg\} (1-r^2)^\alpha r dr\\
&=&2a_0\int_0 ^1 r(1-r^2)^\alpha dr
=f(0)\int_0 ^1 r(1-r)^\alpha dr\\
&=&f(0)\frac{1}{\alpha+1}\ .
\end{eqnarray*}
Multiplying by $\alpha+1$ we get the lemma.
\end{proof}

We now use a disc automorphism to get a representation formula for $f$ at any $z\in\Delta$.

\begin{lemma}\label{lm2.3}
Let $f\in H^\infty(\Delta)$ and $\alpha>-1$. Then
\begin{equation}\label{eq2.7}
f(z)=(\alpha+1)\int_\Delta \frac{f(\zeta)}{(1-\overline{\zeta}z)^{\alpha+2}}(1-|\zeta|^2)^\alpha \ dm(\zeta).
\end{equation}
\end{lemma}
\begin{proof} For $z\in \Delta$ define the disc automorphism $\phi_z :\Delta \longrightarrow \Delta$ as $$\phi_z(w)=\frac{z-w}{1-\overline{z}w}\ .$$
Clearly $\phi_z(0)=z$. Whatsmore, it's easy to establish the properties
\begin{eqnarray*}
&& \phi^{-1} _z =\phi _z, \\
&& \phi_z ^\prime (w) = \frac{1-|z|^2}{(1-\overline{z}w)^2} ,\\
&& 1-|\phi_z(w)|^2=|\phi_z ^\prime(w)|(1-|w|^2).
\end{eqnarray*}
Using these properties and Lemma \ref{lm2.2}, we can write, for any $z\in \Delta$,
\begin{eqnarray*}
f(z)&=&(f\circ\phi_z)(0)=(\alpha+1)\int_\Delta (f\circ\phi_z)(\zeta)(1-|\zeta|^2)^\alpha \ dm(\zeta) \\ 
&=&(\alpha+1)\int_{\phi_z(\Delta)}f(\zeta)|\phi_z ^\prime(\zeta)|^2	(1-|\phi_z(\zeta)|^2)^\alpha \ dm(\zeta) \\
&=&(\alpha+1)\int_\Delta f(\zeta)\frac{(1-|z|^2)^2}{|1-\overline{z}\zeta|^4} \ 
\frac{(1-|z|^2)^\alpha(1-|\zeta|^2)^\alpha}{|1-\zeta \overline{z}|^{2\alpha}} \ dm(\zeta)\\
&=&(\alpha+1)(1-|z|^2)^{\alpha+2}\int_\Delta f(\zeta) \frac{(1-|\zeta|^2)^\alpha}{|1-\zeta\overline{z}|^{2\alpha+4}} \ dm(\zeta) \ .
\end{eqnarray*}
Now, write the above identity for the function $g$, defined as $g(\zeta)=(1-\overline{z}\zeta)^{\alpha+2}f(\zeta)$, in the place of $f$ to get
$$(1-|z|^2)^{\alpha+2}f(z) = (\alpha+1)(1-|z|^2)^{\alpha+2}\int_\Delta \frac{ f(\zeta)(1-|\zeta|^2)^\alpha}{(1-\overline{\zeta}z)^{\alpha+2}} \ dm(\zeta)\ .$$
This proves the lemma.
\end{proof}

So formula (\ref{eq2.7}) defines the family of representations we were seeking for, at least for nice functions in $H^\infty(\Delta)$. The next step is to try to extend this formula to some bigger space that will hopefully incude the spaces $A^p(\Delta)$, for $1\leq p < \infty$. This gives rise to the \textbf{weighted Bergman spaces}.

\begin{definition} For $\alpha>-1$ we define the family of measures $dm_\alpha(z)$ as 
\begin{equation}\label{eq2.8}
dm_\alpha(z) = (\alpha+1)(1-|z|^2)^\alpha dm(z).
\end{equation}
Let $0<p<\infty $. The \textbf{weighted Bergman spaces} $A^p _\alpha (\Delta)$ are defined as 
\begin{eqnarray}\label{eq2.9}
A^p _\alpha (\Delta)=\bigg\{ f\in H(\Delta) : \|f\| _{A_\alpha ^p (\Delta ) }  ^p = \int _{\Delta}|f(z)|^p dm_\alpha(z)<\infty \bigg\}.
\end{eqnarray} 
\end{definition}

The spaces $A^p _\alpha (\Delta)$ consist of the functions that are in $L^p _\alpha (\Delta)=L^p(\Delta,dm_\alpha)$ and are analytic on the unit disc $\Delta$. Thus, the spaces $A^p _\alpha (\Delta)$ are closed subspaces of $L^p _\alpha (\Delta)$ and one can write down a series of results analogous to the ones we have seen for the usual Bergman spaces. Let us first extend formula (2.7) to the case of $A^p _\alpha(\Delta)$ functions. This is the analogous of Proposition \ref{pr2.1} for the weighted case.

\begin{proposition}\label{pr2.4}
Let $\alpha>-1$ and $f\in L^p_\alpha(\Delta)$, $1\leq p < \infty $. Define the function $P_\alpha(f)$ on $\Delta$ as
\begin{equation}\label{eq2.10}
P_\alpha(f)(z)=\int_\Delta \frac{f(\zeta)}{(1-\overline{\zeta}z)^{\alpha+2}}\ dm _\alpha(\zeta), \ \ \ z\in \Delta .
\end{equation}
Then,\\
\noindent(i) The function $P_\alpha(f)$ is a well defined analytic function on $\Delta$.\\
\noindent(ii) If in addition $f\in A^p_\alpha(\Delta)$, $1\leq p <\infty$ and $z\in \Delta$, we have that
\begin{equation}\label{eq2.11}
P_\alpha(f)(z)=f(z).
\end{equation}
\end{proposition} 

The proof is just a repetition of the arguments in the proof of Proposition \ref{pr2.1} and so it's omitted.

It should be clear by now that the operators $P_\alpha$ play the role of the Bergman Projection in the weighted case. One can easily see that $P_\alpha$ is the orthogonal projection of $L^2_\alpha(\Delta)$ onto $A^2_\alpha(\Delta)$ and that the operator $P_\alpha:L^p_\alpha(\Delta)\longrightarrow A^p_\alpha(\Delta)$ is bounded whenever $1<p<\infty$. Remember however that we seek to find a bounded projection of $L^1(\Delta)$ onto $A^1(\Delta)$ in the \textbf{non-weighted} case so this is not very helpful.

When $\alpha \geq 0$, the weighted Bergman spaces contain the usual Bergman spaces, that is, $A^p(\Delta)\subseteq A^p_\alpha(\Delta)$, and of course $L^p(\Delta)\subseteq L^p_\alpha(\Delta)$ while the non-weighted case corresponds to the value $\alpha=0$. This means that formula (\ref{eq2.10}) still holds if $\alpha>0$ and $f\in L^1(\Delta)$. This is more promising! What we really need is a thorough description of when the operator $P_\alpha$ is a bounded projection of $L^p(\Delta)$ onto $A^p(\Delta)$. We already know this description when $\alpha=0$ (this is Theorem \ref{thm2.1}). The following Theorem gives an answer in the general case $\alpha>-1$.

\begin{theorem}\label{thm2.2} Let $\alpha>-1$ and $1\leq p <\infty$. For $f\in L^p(\Delta)$ we define the function $P_\alpha(f)$ as
$$P_\alpha(f)=(\alpha+1)\int_\Delta\frac{f(\zeta)}{(1-\overline{\zeta}z)^{\alpha+2}}
(1-|\zeta|^2)^\alpha \ dm(\zeta) .$$
Then $P_\alpha$ is a bounded operator from $L^p(\Delta)$ onto $A^p(\Delta)$ if and only if $p(\alpha+1)>1$.
\end{theorem}
\noindent\textbf{Remarks.}
\noindent(i) When $\alpha=0$, this is just  a  repetition of the statement of Theorem \ref{thm2.1}, that is that the Bergman Projection is bounded on $L^p(\Delta)$ if and only if $1<p<\infty$.  \smallskip \\
\noindent(ii) The theorem tells us that if $\alpha>0$ then $P_\alpha$ carries boundedly $L^p(\Delta)$ onto $A^p(\Delta)$ for any $1\leq p < \infty$. More specifically this means that there exists a bounded projection of $L^1(\Delta)$ onto $A^1(\Delta)$ which is what we were seeking for. Remember that this is not the case for the Hardy space $H^1(\Delta)$. \smallskip \\
\noindent(iii) There exists an even more general version of this theorem that says that 
if $-1<\alpha,\beta<\infty$ then the operator $P_\alpha : L^p_\beta(\Delta)\longrightarrow A^p_\beta(\Delta)$ is a bounded projection of $L^p_\beta(\Delta)$ onto $A^p_\beta(\Delta)$ if and only if $\beta+1<(\alpha+1)p$. Since we wont need the result in this generality, we will only give the proof for the case $\beta=0$.

\begin{proof}{Case $p=1$}. Let us first show that if $P_\alpha$ is bounded on $L^1(\Delta)$ if and only if $\alpha>0$. To that end we will use the fact that the adjoint operator of $P_\alpha$, $P^* _\alpha$, is bounded on $L^\infty(\Delta)$ if and only if $P_\alpha$ is bounded on $L^1(\Delta)$.
The operator $P_\alpha ^*$ is defined by means of the "inner product"
$$\langle f,g \rangle =\int_\Delta f(z)\overline{g(z)} dm(z)$$	
where $f\in L^1(\Delta)$ and $g\in L^\infty(\Delta)$. That is, we define the operator $P_\alpha ^*$ on $L^\infty(\Delta)$ so that for every pair of functions $f\in L^1(\Delta)$ and $g\in L^\infty(\Delta)$ we have that
$$\langle P_{\alpha}(f),g \rangle =\langle f,P^* _{\alpha}(g)\rangle .$$
The left hand-side inner product can be calculated as follows.
\begin{eqnarray*}
\langle P_{\alpha}(f),g \rangle &=&\int_\Delta P_\alpha(f)(\zeta)\overline{g(\zeta)}dm(\zeta)=  	
\int_\Delta \int_\Delta \frac{f(z)}{(1-\overline{z}\zeta)^{2+\alpha}}\ dm_\alpha(z) \overline{g(\zeta)}\ dm(\zeta) \\ \\
&=&\int_\Delta f(z)\int_\Delta \frac{\overline{g(\zeta)}}{(1-\overline{z}\zeta)^{2+\alpha}} \ dm(\zeta) \ dm_\alpha(z)\\
\\ &=&\int_\Delta f(z)\overline{\int_\Delta \frac{(\alpha+1)g(\zeta)(1-|z|^2)^\alpha}{(1-\overline{\zeta}z)^{2+\alpha}} \ dm(\zeta) }\ dm_\alpha(z).
\end{eqnarray*}
We have thus found an explicit formula for the operator $P_\alpha ^*$,
\begin{equation}\label{eq2.12}
P^* _\alpha (g)(z)=\int_\Delta \frac{(\alpha+1)g(\zeta)(1-|z|^2)^\alpha}{(1-\overline{\zeta}z)^{2+\alpha}} \ dm(\zeta). 
\end{equation}
It is now easy to see that $P_\alpha ^*$ is bounded on $L^\infty(\Delta)$ if and only if \begin{equation}\label{eq2.13}
\sup_{z\in\Delta} (1-|z|^2)^\alpha \int _\Delta \frac{dm(\zeta)}{|1-\overline{\zeta}z|^{2+\alpha}}<\infty 
\end{equation}
(just test $P^* _\alpha$ against the function $g(\zeta)=\frac{(1-z\overline{\zeta})^{2+\alpha}}{|1-\overline{\zeta} z|^{2+\alpha}}$).

Let us next show that equation (\ref{eq2.13}) holds if and only if $\alpha>0$. Indeed, suppose that $\alpha>0$. Then, 
\begin{eqnarray*}
(1-|z|^2)^\alpha \int _\Delta \frac{dm(\zeta)}{|1-\overline{\zeta}z|^{2+\alpha}}&=&
(1-|z|^2)^\alpha \int _0 ^1 \frac{1}{\pi}\int_0 ^{2\pi} \frac{1}{|1-\overline{z}re^{i\theta}|^{2+\alpha}} d\theta \ r dr \\
&=&(1-|z|^2)^\alpha \sum _{n=0} ^\infty \Big(\frac{\Gamma(n+\frac{\alpha}{2}+1)}{n!\Gamma(\frac{\alpha}{2}+1)}\Big)^2\int_0 ^1 r^{2n+1}dr |z|^{2n} \\
&=&\frac{(1-|z|^2)^\alpha}{2} \sum _{n=0} ^\infty \Big(\frac{\Gamma(n+\frac{\alpha}{2}+1)}{n!\Gamma(\frac{\alpha}{2}+1)}\Big)^2\frac{1}{n+1}|z|^{2n} \\
&\leq& C \frac{(1-|z|^2)^\alpha}{2\Big(\Gamma(\frac{\alpha}{2}+1)\Big)^2} \sum _{n=0} ^\infty (n+1)^{\alpha-1}|z|^{2n}\leq C^\prime ,
\end{eqnarray*}
where $C, C^\prime$ are absolute constants. On the other hand, if $\alpha=0$, we have that
$$\int_\Delta\frac{1}{|1-\overline{\zeta}z|^2} dm(\zeta)=\sum_{n=0} ^\infty 2 \int_0 ^1 r^{2n+1} dr |z|^{2n}=\sum_{n=0} ^\infty \frac{1}{n+1}|z|^{2n}\geq c \log\frac{1}{1-|z|},$$ for some absolute constant $c$. Finally, if $-1<\alpha<0$ then
$$(1-|z|^2)^\alpha	 \int _\Delta \frac{dm(\zeta)}{|1-\overline{\zeta}z|^{2+\alpha}}
\geq (1-|z|^2)^\alpha \int _\Delta \frac{dm(\zeta)}{(1+|z|)^{2+\alpha}}\geq \frac{1}{2^{2+\alpha}}(1-|z|^2)^\alpha .$$ This shows that if $P_\alpha$ is bounded on $L^1(\Delta)$ then we must have that $\alpha>0$.

Now let $\alpha>0$. We will show that $P_\alpha$ is bounded on $L^1(\Delta)$. We have that
\begin{eqnarray*}
\int_\Delta |P_\alpha(f)(z)|dm(z)&=&\int_\Delta \bigg|(\alpha+1)\int_\Delta \frac{f(\zeta)(1-|\zeta|^2)^\alpha}{(1-\overline{\zeta}z)^{2+\alpha}}\ dm(\zeta)\bigg|\  dm(z) \\ 
&\leq& (\alpha+1)\int_\Delta \int_\Delta \frac{|f(\zeta)|}{|1-\overline{\zeta}z|^{2+\alpha}}(1-|\zeta|^2)^\alpha dm(\zeta)dm(z) \\
&=& (\alpha+1)\int_\Delta|f(\zeta)| \int_\Delta \frac{1}{|1-\overline{\zeta}z|^{2+\alpha}}dm(z)(1-|\zeta|^2)^\alpha dm(\zeta) .
\end{eqnarray*}
However, 
\begin{eqnarray*}
\int_\Delta \frac{1}{|1-\overline{\zeta}z|^{2+\alpha}}dm(z)&=&\int_\Delta \frac{1}{|(1-\overline{\zeta}z)^{\frac{2+\alpha}{2}}|^2}dm(z)\\
&=& \frac{1}{\big(\Gamma(\frac{\alpha}{2}+1)\big)^2}\sum_{n=0} ^\infty \Bigg(\frac{\Gamma(n+\frac{\alpha}{2}+1)}{n!} \Bigg)^2\frac{1}{n+1}|\zeta|^{2n}
\\ &\leq& c\ \frac{\Gamma(\alpha)}{\Gamma(\frac{\alpha}{2}+1)^2} \sum_{n=0} ^\infty \frac{\Gamma(n+\alpha)}{n!\Gamma(\alpha)}|\zeta|^{2n}=\frac{\Gamma(\alpha)}{\Gamma(\frac{\alpha}{2}+1)}\frac{1}{(1-|\zeta|^2)^\alpha}.
\end{eqnarray*}
Hence,
$$\int_\Delta |P_\alpha(f)(z)|dm(z)\leq C\ \|f\|_{L^1(\Delta)},$$
for some constant $C$ depending only on $\alpha$. This finishes the case $p=1$.

\smallskip
\noindent {Case $1<p<\infty$}. Suppose now that $P_\alpha \in B(L^p(\Delta,dm))$. This implies that $P^* _\alpha \in B(L^q(\Delta,dm))$ where $q=\frac{p}{p-1}$ is the conjugate exponent of $p$. We will show that $p(\alpha+1)>1$. Suppose, for the sake of contradiction, that $p(\alpha+1)\leq 1$ which is equivalent to $\alpha \leq -\frac{1}{q}$.

Remember that the operator $P_\alpha ^*$ is described by equation (\ref{eq2.12}). Taking $g(\zeta)=1\in L^q(\Delta,dm)$ we have
\begin{eqnarray*}
P_\alpha ^*(g)(z)&=&\int_\Delta \frac{(\alpha+1)(1-|z|^2)^\alpha}{(1-\overline{\zeta}z)^{2+\alpha}}dm(\zeta)
\\&=&(\alpha+1)(1-|z|^2)^\alpha \int_0 ^1\frac{1}{\pi}\int_0 ^{2\pi} \frac{1}{(1-zre^{-i\theta})^{2+\alpha}}\ d\theta \ rdr\\
&=&(\alpha+1)(1-|z|^2)^\alpha \int_0 ^1 
\frac{1}{\pi}\int_0 ^{2\pi}\sum_{n=0} ^\infty \frac{\Gamma(n+\alpha+2)}{n!\Gamma(\alpha+2)} r^n z^n e^{-in\theta }d\theta rdr\\
&=&(\alpha+1)(1-|z|^2)^\alpha \int_0 ^1\sum_{n=0} ^\infty 
\frac{1}{\pi}\int_0 ^{2\pi} e^{-in\theta }d\theta r^n z^n rdr\\
&=&(\alpha+1)(1-|z|^2)^\alpha \int_0 ^1 rdr =(\alpha+1)(1-|z|^2)^\alpha .
\end{eqnarray*}
However, the function $(\alpha+1)(1-|z|^2)^\alpha$ is not in $L^q(\Delta,dm)$ and hence we have a contradiction. Indeed, if $\alpha<-\frac{1}{q}$ then 
\begin{eqnarray*}
(\alpha+1)^q\int_\Delta(1-|z|^2)^{\alpha q}dm(z)=(\alpha+1)^q\int_0 ^1 x^{\alpha q }dx\\
=(\alpha+1)^{\alpha q}\bigg[\frac{1}{\alpha q+1} \big(1-\lim_{x\rightarrow 0} \frac{1}{x^{-1-\alpha q}} \big)\bigg]=\infty.
\end{eqnarray*}
On the other hand, if $\alpha=-\frac{1}{q}$ then
$$(\alpha+1)^q\int_\Delta(1-|z|^2)^{-1}dm(z)=(\alpha+1)^q\int_0 ^1 \frac{1}{x}dx =\infty.$$
In order to complete the proof we shall need a boundedness criterion known as \emph{Schur's Test}. 
\begin{theorem}(Schur's Test)\label{thm2.3}
Let $(\mathbb{X},\mu)$ be a measure space and $K:\mathbb{X}\times \mathbb{X}\rightarrow \R^+ $ a nonnegative measurable function. For $1<p<\infty$ and $f\in L^p(\mathbb{X},d\mu) $ we define
$$T(f)(x)=\int_{\mathbb{X}}K(x,y)f(y) d\mu(y), \ x\in \mathbb{X}.$$
Suppose that there exists a positive constant $C>0$ and a positive measurable function $h$ on $\mathbb{X}$ such that
\item(a)For $\mu-$almost every $x\in\mathbb{X}$,$$ \int _{\mathbb{X}}K(x,y) h(y)^q d\mu (y)\leq C h(x)^q \ .$$
\item(b)For $\mu-$almost every $y\in\mathbb{X}$,$$ \int _{\mathbb{X}}K(x,y) h(x)^p d\mu (x)\leq C h(y)^p \ .$$
Then, $T\in B(L^p(\mathbb{X},d\mu))$ and $\|T\|\leq C$.
\end{theorem}
\noindent\textbf{Remark.} Suppose that $(\mathbb{X},\mu)$ is a measure space and $K:\mathbb{X}\times \mathbb{X}\rightarrow \R^+ $ a nonnegative measurable function. Suppose further that there exists some positive constant $C$ and some positive measurable function $h$, defined on $\mathbb{X}$, such that
\item(i)$\int_\mathbb{X} K(x,y)h(y)d\mu(y)\leq C h(x)$ for $\mu$-almost every $x\in\mathbb{X}$.
\item(ii)$\int_\mathbb{X} K(x,y)h(x)d\mu(x)\leq C h(y)$ for $\mu$-almost every $y\in\mathbb{X}$.
\smallskip

\noindent Then, the operator $T(f)(x)=\int_\mathbb{X} K(x,y)f(y)d\mu(y)\in B(L^2(\mathbb{X},d\mu))$ with $\|T\|\leq C$.

\medskip
We postpone the proof of Schur's test until after the end of the proof of Theorem \ref{thm2.2}

\medskip

Let's fix some $p\in(1,\infty)$ and set $q=\frac{p}{p-1}$. Suppose that $p(\alpha+1)>1$. We will now use Schur's test in order to show that $P_\alpha$ is a bounded operator from $L^p(\Delta)$ to $A^p(\Delta)$.

We set $d\mu(\zeta)=(1-|\zeta|^2)^\alpha dm(\zeta)$, $h(z)=(1-|z|^2)^{-\frac{1}{pq}}$ and $K(z,\zeta)=\frac{1}{|1-\overline{\zeta}z|^{2+\alpha}} \ .$ We have that
\begin{eqnarray*}
\int_\Delta K(z,\zeta)h(\zeta)^q d\mu(\zeta)&=&
\int_\Delta \frac{1}{|1-\overline{\zeta}z|^{2+\alpha}}h(\zeta)^q(1-|\zeta|^2)^{\alpha}dm(\zeta)\\
&=&\int_\Delta \frac{(1-|\zeta|^2)^{\alpha-\frac{1}{p}}}{|1-\overline{\zeta}z|^{2+\alpha}}d\mu(\zeta).
\end{eqnarray*}
However, since
$\frac{1}{(1-\zeta\overline{z})^\frac{2+\alpha}{2}}=\sum_{n=0}^\infty\frac{\Gamma(n+\frac{\alpha}{2}+1)}{n!\Gamma(\frac{\alpha}{2}+1)}\overline{z}^n\zeta^n ,$ we get that, 
\begin{eqnarray*}
\int_\Delta \frac{(1-|\zeta|^2)^{\alpha-\frac{1}{p}}}{|1-\overline{\zeta}z|^{\alpha+2}}dm(\zeta)
&\leq & \sum _{n=0}^\infty \bigg(\frac{\Gamma(n+\frac{\alpha}{2}+1)}{n!\Gamma(\frac{\alpha}{2}+1)} \bigg)^2\int _0 ^1 r^n (1-r)^{\alpha-\frac{1}{p}} dr |z|^{2n} \\
&=& \sum _{n=0}^\infty \bigg(\frac{\Gamma(n+\frac{\alpha}{2}+1)}{n!\Gamma(\frac{\alpha}{2}+1)} \bigg)^2 \frac{\Gamma(n+1)\Gamma(\alpha+1-\frac{1}{p})}{\Gamma(n+\alpha+2-\frac{1}{p})}|z|^{2n}\\
&=&\frac{\Gamma(\alpha+\frac{1}{q})}{\Gamma(\frac{\alpha}{2}+1)^2}
\sum_{n=0}^\infty \bigg(\frac{\Gamma(n+\frac{\alpha}{2}+1)}{n!} \bigg)^2 
\frac{1}{\frac{\Gamma(n+\alpha+1+\frac{1}{q})}{n!}}|z|^{2n}\\
&\leq& C \frac{\Gamma(\alpha+\frac{1}{q})}{\Gamma(\frac{\alpha}{2}+1)^2} 
\sum_{n=0}^\infty \frac{\Gamma(n+1-\frac{1}{q})}{n!}|z|^{2n}=
C \frac{\Gamma(\alpha+\frac{1}{q})}{\Gamma(\frac{\alpha}{2}+1)^2}\frac{1}{(1-|z|)^\frac{1}{p}}\\
&=&C_\alpha h(z)^q,
\end{eqnarray*}
for some constant $C_\alpha$, depending only on $\alpha$.

Using a similar calculation, we show that
$$\int_\Delta K(z,\zeta)h(z)^p d\mu(z)\leq C_\alpha h(\zeta)^p .$$
Now, Schur's test tells us that $P_\alpha\in B(L^p(\Delta),A^p(\Delta))$ and that $\|P_\alpha\|\leq C_\alpha$.
\end{proof}
\noindent\textit{Proof of theorem {\ref{thm2.3}}}. Let $f\in L^p(\mathbb{X},d\mu)$. Then,

$$|T(f)(x)|\leq  \int _\mathbb{X}K(x,y)h(y)\frac{1}{h(y)}|f(y)|d\mu(y).$$
H{\"o}lder's inequality now yields
\begin{eqnarray}
 |T(f)(x)|&\leq& \bigg( \int _\mathbb{X} K(x,y)h(y)^q d\mu(y)\bigg)^\frac{1}{q}
  \bigg( \int _\mathbb{X} K(x,y)\frac{|f(y)|^p}{h(y)^p} d\mu(y)\bigg)^\frac{1}{p}\nonumber \\
 &\leq& C^\frac{1}{q} h(x) \bigg( \int _\mathbb{X} K(x,y)h(y)^{-p}|f(y)|^p d\mu(y)\bigg)^\frac{1}{p} \label{eq2.14}
\end{eqnarray}
for $\mu$-almost every $x\in\mathbb{X}$. Integrating with respect to $x$ and using Fubini we get
\begin{eqnarray*}
 \int_{\mathbb{X}}|T(f)(x)|^pd\mu(x) &\leq& C^\frac{p}{q}\int_\mathbb{X} h(x)^p \int_\mathbb{X} K(x,y)h(y)^{-p}|f(y)|^pd\mu(y)\ d\mu(x) \\
 &=& C^\frac{p}{q}\int_\mathbb{X}\int_\mathbb{X}K(x,y)h(x)^pd\mu(x)h(y)^{-p}|f(y)|^p d\mu(y)\\
 &\leq & C^\frac{p}{q}+1\int_\mathbb{X}|f(y)|^p d\mu(y),
\end{eqnarray*}
where the last inequality follows by (\ref{eq2.14}). Thus we have that
$$ \|T(f)\|^p _{L^p(\mathbb{X}d\mu)}\leq C^{\frac{p}{q}+1}\|f\|^p _{L^p(\mathbb{X}d\mu)}.$$
Taking $p$-th roots copletes the proof.\qed
\medskip


\section {A characterization of $A^p$ in terms of derivatives}
\medskip

Suppose that $f$ is an analytic function and fix some $p\in[0,\infty)$. Let $n$ be a positve integer. We want to find a condition for $f^{(n)}$ that assures that $f\in A^p(\Delta)$. We actually get something better, that is, a characterization of the space $A^p(\Delta)$ in terms of derivatives. Our main result for this section is the following.

\begin{theorem}\label{thm2.4}
Let $n$ be a positive integer greater than 1 and $1\leq p<\infty $. Suppose that $f\in H(\Delta)$. Then $f\in A^p(\Delta)$ if and only if $(1-|z|^2)f^{(n)}(z)\in L^p(\Delta)$.
\end{theorem}

The proof will be done in several steps. Let's begin with the necessity of the condition of Theorem \ref{thm2.4}.

\begin{lemma}\label{lm2.4}Let $n$ be a positive integer greater than 1 and $1\leq p<\infty $. Suppose that $f\in A^p(\Delta)$. Then $(1-|z|^2)^n f^{(n)}\in L^p(\Delta)$.
\end{lemma}
\begin{proof}{Case $p=1$.} Suppose that $f\in A^1(\Delta)$. For $\alpha=1$, equation (\ref{eq2.10}) of Proposition \ref{pr2.4} emplies 
$$ f(z)=P_\alpha(f)(z)=\int_\Delta \frac{f(\zeta)}{(1-\overline{\zeta}z)^{\alpha+2}}\ dm _\alpha(\zeta)=2\int_\Delta \frac{(1-|\zeta|^2)}{(1-\overline{\zeta}z)^3}f(\zeta)dm(\zeta) .$$ 
Differentiating $n$ times we get
$$ f^{(n)}(z)=(n+2)!\int_\Delta \frac{(\overline{\zeta})^n(1-|\zeta|^2)}{(1-\overline{\zeta}z)^{n+3}}f(\zeta)dm(\zeta).$$
Thus,
\begin{eqnarray*}
\int_\Delta(1-|z|^2)^n|f^{(n)}(z)|dm(z)&\leq & (n+2)!\int_\Delta (1-|z|^2)^n
\int_\Delta \frac{(1-|\zeta|^2)}{|1-\overline{\zeta}z|^{n+3}}|f(\zeta)|dm(\zeta)dm(z)\\
 &\leq& (n+2)!\int_\Delta (1-|\zeta|^2)|f(\zeta)|
\int_\Delta \frac{(1-|z|^2)^n}{|1-\overline{\zeta}z|^{n+3}}dm(z) dm(\zeta)\\
&\leq& C \int_\Delta (1-|\zeta|^2)|f(\zeta)|\frac{1}{(1-|\zeta|)}dm(\zeta)\\&=&C\|f\|_{A^1(\Delta)}.
\end{eqnarray*}
\smallskip
\noindent{Case $1<p<\infty$.} Suppose that $f\in A^p(\Delta)$. By Proposition \ref{pr2.1}
we have that
$$f(z)=\int_\Delta \frac{f(\zeta)}{(1-\overline{\zeta}z)^2}dm(\zeta).$$
Differentiating $n$ times we get
\begin{eqnarray*}
(1-|z|^2)^nf^{(n)}(z)=(n+1)!(1-|z|^2)^n\int _\Delta \frac{(\overline{\zeta})^nf(\zeta)}{(1-\overline{\zeta}z)^{n+2}}dm(\zeta)=n!(P_n ^*\circ S_n)(f)(z),
\end{eqnarray*}
where $S^n(f)(z)={\overline{z}}^nf(z)$ and $P_n^*$ is described by equation (\ref{eq2.12}). Consequently
$$\int_\Delta |1-|z|^2|^{pn}|f^{(n)}(z)|^pdm(z)\leq n!\|P_n ^*\| \ \|S_n\| \ \|f\|_{A^p(\Delta)}\leq C_n \|f\|_{A^p(\Delta)}.$$
Remember that $P_n ^*\in B(L^p(\Delta))$ by Theorem \ref{thm2.2}.
\end{proof}

To prove the sufficiency in Theorem \ref{thm2.4} we need the following lemma.

\begin{lemma}\label{lm2.5} Let $f\in H(\Delta)$ and $n$ be a positive integer greater than 1, such that 
\item(i) $(1-|z|^2)^n f^{(n)}(z)\in L^1(\Delta)$, and
\item(ii)$f(0)=f'(0)=f''(0)=\cdots=f^{(2n-1)}(0)=0.$

\smallskip
\noindent Then, for every $z\in \Delta$, 

$$f(z)=\frac{1}{n!}\int_\Delta \frac{(1-|\zeta|^2)^nf^{(n)}(\zeta)}{(\overline{\zeta})^n(1-\overline{\zeta}z)^2}dm(\zeta).$$
\end{lemma}

\begin{proof} First notice that condition (ii) implies
\begin{eqnarray*}
f^{(n)}(z)&=&\sum_{m=2n} ^\infty m(m-1)\ldots(m-n+1)a_m z^{m-n}\\
&=&z^n\sum_{m=2n} ^\infty m(m-1)\ldots(m-n+1)a_m z^{m-2n}.
\end{eqnarray*}
As a result, $\frac{f^{(n)}(z)}{z^n}\in H(\Delta)$.

For $z\in\Delta$, we write
\begin{eqnarray*}
\int_\Delta\frac{(1-|\zeta|^2)|f^{(n)}(\zeta)|}{|1-\overline{\zeta}z|^2|\overline{\zeta}^n|}dm(\zeta)&\leq& \frac{1}{(1-|z|)^2}\int_\Delta (1-|\zeta|^2)^n\bigg|\frac{f^{(n)}(\zeta)}{\overline{\zeta}^n}\bigg|dm(\zeta)\\
&\leq&\frac{1}{(1-|z|)^2}\bigg\{\int_{\overline{\Delta_{\frac{1}{2}}}}(1-|\zeta|^2)^n \bigg|\frac{f^{(n)}(\zeta)}{\zeta^n}\bigg|dm(\zeta)\\&+&\int_{\Delta  -\overline{\Delta_{\frac{1}{2}}}}(1-|\zeta|^2)^n \bigg|\frac{f^{(n)}(\zeta)}{\zeta^n}\bigg|dm(\zeta)\bigg\}.
\end{eqnarray*}
For the first integral 	notice that $\sup_{\zeta \in \overline{\Delta_\frac{1}{2}}}\bigg|\frac{f^{(n)}(\zeta)}{\zeta^n}\bigg|<\infty.$ Hence $$\int_{\overline {\Delta_\frac{1}{2}}} (1-|\zeta|^2)^n|f^{(n)}(\zeta)|dm(\zeta)<\infty . $$
On the other hand
$$\int_{\Delta  -\overline{\Delta_{\frac{1}{2}}}}(1-|\zeta|^2)^n \bigg|\frac{f^{(n)}(\zeta)}{\zeta^n}\bigg|dm(\zeta)\leq2^n \int_{\Delta  -\overline{\Delta_{\frac{1}{2}}}}(1-|\zeta|^2)^n|f^{(n)}(\zeta)|dm(\zeta)<\infty$$
due to hypothesis (i). This shows that the integral 
$\int_\Delta \frac{(1-|\zeta|^2)^nf^{(n)}(\zeta)}{(\overline{\zeta})^n(1-\overline{\zeta}z)^2}dm(\zeta)$ exists. What we actually showed is that the function $F(z)=\frac{(1-|z|^2)^nf^{(n)}(z)}{\overline{z}^n}$ is in $L^1(\Delta).$ But this means that the integral
$$g(z)=\int_\Delta\frac{F(\zeta)}{(1-\overline{\zeta}z)^2}dm(\zeta)=\frac{1}{n!}\int_\Delta \frac{(1-|\zeta|^2)^nf^{(n)}(\zeta)}{(1-\overline{\zeta}z)^2\overline{\zeta}^n}dm(\zeta) $$
defines an analytic function on $\Delta$. Hence
\begin{eqnarray*}
g^{(n)}(z)&=&(n+1)\int_\Delta \frac{(1-|\zeta|^2)^nf^{(n)}(\zeta)}{(1-\overline{\zeta}z)^{n+2}\overline{\zeta}^n}\overline{\zeta}^ndm(\zeta)\\&=&(n+1)\int_\Delta \frac{(1-|\zeta|^2)^nf^{(n)}(\zeta)}{(1-\overline{\zeta}z)^{n+2}}dm(\zeta)=P_n (f^{(n)})(z).
\end{eqnarray*}
Since $f^{(n)}\in A^1 _n(\Delta)$, part (ii) of Proposition \ref{pr2.4} gives that $P_n(f^{(n)})(z)=f^{(n)}(z)$. This means that $f^{(n)}(z)=g^{(n)}(z)$. However, for $0\leq k \leq n-1$, we have that
\begin{eqnarray*}
g^{(k)}(z)=\frac{(k+1)!}{n!}\int_\Delta \frac{f^{(n)(\zeta)}(1-|\zeta|^2)\overline{\zeta}^k}{(1-\overline{\zeta}z)^{n+2}\overline{\zeta}^n}dm(\zeta).
\end{eqnarray*}
But this implies
\begin{eqnarray*}
g^{(k)}(0)&=&\frac{(k+1)!}{n!}\int_\Delta \frac{f^{(n)(\zeta)}(1-|\zeta|^2)}{\overline{\zeta}^{n-k}}dm(\zeta)\\
&=&\frac{(k+1)!}{n!}\int	_0 ^1 \sum_{m=2n} ^\infty m(m-1)\ldots(m-n+1)r^{m+k-2n}\frac{1}{2\pi}\int_0 ^{2\pi}e^{i(m-k)\theta}d\theta \ r dr	\\&=&0=f^{(k)}(0).
\end{eqnarray*}
Hence $f=g$ and the proof is complete.\end{proof}

We are now ready to prove Theorem \ref{thm2.4}
\medskip

\noindent\textit{Proof of Theorem \ref{thm2.4}.} Let $1\leq p <\infty$ and $n$ be a postive integer, greater than 1. Suppose that $f\in H(\Delta)$ is such that $(1-|z|^2)^nf^{(n)}(z)\in L^p(\Delta)$. Suppose further that
$$f(0)=f'(0)=f''(0)=\cdots=f^{(2n)}(0)=0.$$
As we've seen in Lemma \ref{lm2.5}, if we set $F(z)=\frac{(1-|z|^2)^nf^{(n)}(z)}{\overline{z}^n}$ is in $L^1(\Delta)$, then 
$$f(z)=\frac{1}{n!}\int_\Delta \frac{(1-|\zeta|^2)^n f^{(n)}(\zeta)}{(1-\overline{\zeta}z)^2\overline{\zeta}^n}dm(\zeta)=P(f)(z).$$

If $1<p<\infty$, then $P\in B(L^p(\Delta,dm))$ and hence $f\in A^p(\Delta)$.

Let $p=1$. Then, for $\alpha>0$, the operator
$$P_\alpha(f)(z)=(\alpha+1)\int_\Delta\frac{f(\zeta)}{(1-\overline{\zeta}z)^{\alpha+2}}(1-|\zeta|^2)dm(\zeta)$$
is bounded on $L^1(\Delta)$.\qed

\newpage
\newpage

\end{document}